\documentclass[12pt,a4paper]{article}

\usepackage[margin=1in]{geometry}
\usepackage{graphicx}    % For including figures
\usepackage{caption}     % Better control over captions
\usepackage{enumitem}    % Custom bullet lists
\usepackage{setspace}    % For spacing
\usepackage{titlesec}    % For section formatting
\usepackage{float}
\usepackage{amsmath,amssymb,mathtools,amsthm}
\mathtoolsset{showonlyrefs}
\usepackage{xcolor}
\usepackage{authblk}
\usepackage{dsfont}
\usepackage[colorlinks,
            pdffitwindow=false,
            plainpages=false,
            pdfpagelabels=true,
            pdfpagemode=UseOutlines,
            pdfpagelayout=SinglePage,
            bookmarks=false,
            colorlinks=true,
            hyperfootnotes=false,
            linkcolor=blue,
            urlcolor=blue!50!black,
            citecolor=green!50!black]
{hyperref}
\theoremstyle{definition}
\newtheorem{remark}{Remark}

\newcommand{\R}{\mathds{R}}
\newcommand{\N}{\mathds{N}}

\newcommand{\E}{\mathds{E}}
\newcommand{\norm}[1]{\left\lVert#1\right\rVert}

\newcommand{\1}{\mathds{1}}
\renewcommand{\P}{\mathds{P}}

\newcommand{\Id}{\textup{Id}}

\newenvironment{nalign}{
    \begin{equation}
    \begin{aligned}
}{
    \end{aligned}
    \end{equation}
    \ignorespacesafterend
}

\setlist[itemize]{leftmargin=1.5em, label=--}
\titleformat{\section}{\large\bfseries}{\thesection.}{0.5em}{}
\title{Noise-induced enhancement of regime lifetimes -- A data-driven approach using deterministic trajectories}

\begin{document}

\author[1]{Henry Schoeller} %% correspondence author
\author[2]{Robin Chemnitz}
\author[3]{P\' eter Koltai}
\author[4,2]{Maximilian Engel}
\author[1]{Stephan Pfahl}

\affil[1]{Institute of Meteorology, Freie Universität Berlin}
\affil[2]{Institute of Mathematics, Freie Universität Berlin}
\affil[3]{Faculty of Mathematics, Universität Bayreuth}
\affil[4]{KdV Institute of Mathematics, University of Amsterdam}

\maketitle

\begin{abstract}
    We investigate the lifetime of dynamical regimes under the impact of noise motivated by low-dimensional models of the atmosphere. One may expect that the inclusion of noise tends to make the system leave prescribed regions of the state space faster. However, for relevant systems with complexities ranging from phenomenological toy models to reduced models of atmospheric dynamics, this intuition has proven misleading. As long as the noise is sufficiently small, the noisy system stays in regimes of interest on average longer than its deterministic counterpart, an effect we call ``stochastic inertia''. This phenomenon has been observed through extensive numerical simulations for different noise levels. We propose a numerical technique for testing the occurrence of stochastic inertia, constructing, for any fixed noise level, a Markov chain on the set of points given by a  sufficiently long trajectory of the system without noise. The method is shown to correctly predict the presence of stochastic inertia in simple systems, and its utility is demonstrated on a paradigm model of atmospheric dynamics.
\end{abstract}

\tableofcontents

\section{Introduction}

The typical long-term behavior of dissipative dynamical systems is governed by the properties of their lower-dimensional attractor. In some cases it is insightful to partition a system's attractor into a finite number of sets, whose properties differ substantially. Such a partitioning and the system's movement between the sets is often characterized by a separation of (time) scales~\cite{dellnitz1999approximation,SchSa13,BoPaNo14}. In some scientific fields, these different sets are termed ``regimes'', most notably in atmospheric sciences~\cite{hannachiLowFrequencyNonlinearity2017}.

The study of ``weather regimes'', loosely defined as persistent flow configurations usually in the mid-latitudes, is long-standing and well-established. In the tradition of Lorenz~\cite{lorenzDeterministicNonperiodicFlow1963}, observed reoccurrence of a finite set of weather regimes has been associated with the structure of the attractor for the atmospheric system \cite{ghilWavesVsParticles2002, strommenTopologicalPerspectiveWeather2023}. Given the near-infinite dimensionality of the atmosphere, studying its attractor is, of course, intractable, so mathematical models are devised as abstractions~\cite{derembleFixedPointsStable2009, ghilReviewArticleDynamical2023, crommelinMechanismAtmosphericRegime2004}. Atmospheric models of any complexity have been shown to reproduce some regime behavior comparable to the one observed in the real atmosphere, and the simplest ones allow for studying the system's attractor itself \cite{maherModelHierarchiesUnderstanding2019}.

The inclusion of some form of randomness in a dynamical system is a natural endeavor, given that stochasticity is a property central to numerous phenomena in science. In atmospheric models in particular, stochasticity has been introduced to alleviate shortcomings in the representation of subgrid scale processes~\cite{palmerStochasticWeatherClimate2019}. Considering that self-similarity principles and scaling cascades are a property of the turbulent Navier--Stokes equations, the parameterization of such processes based purely on the large scale has been shown to introduce biases. Among the biases mitigated by stochastic parameterizations are those that may be regarded as biases in regime representation in the broadest sense~\cite{weisheimerAddressingModelError2014, subramanianImpactStochasticPhysics2017, christensenStochasticParameterizationNino2017, strommenImpactStochasticParametrisations2018, strommenProgressProbabilisticEarth2019} and in the context of weather regimes in the mid-latitudes in the typical sense~\cite{christensenSimulatingWeatherRegimes2015, dawsonSimulatingWeatherRegimes2015, deinhardProcessorientedUnderstandingImpact2024}. In particular, several studies have shown that the persistence of so-called ``blocking'' regimes is altered upon the inclusion of a  stochastic parametrization scheme. Though some studies report an improvement in the representation~\cite{daviniRepresentationWinterNorthern2021, maddisonImpactModelUpgrades2020}, in some cases, the changes may be obfuscated by compensating effects necessitating retuning of model parameters~\cite{filippucciImpactStochasticPhysics2024, hartungResolutionPhysicsAtmosphere2017}.

Between the justification of stochasticity in the parameterization of subgrid scale processes by first principles and the empirical evidence of improved regime representation lies a knowledge gap, pertaining to the precise physical and mathematical mechanisms involved. A step towards closing this gap has been taken by Dorrington and Palmer \cite{dorrington2023interaction}, who have investigated the increase in regime lifetimes upon the inclusion of additive white noise in a low-order paradigm model of mid-latitude atmospheric flow.  This model was first introduced by Charney and deVore~\cite{charneyMultipleFlowEquilibria1979} (CdV system hereafter), and its configuration used both here and by \cite{dorrington2023interaction} is the same as in~\cite{crommelinMechanismAtmosphericRegime2004}. More specifically, for a regime of the CdV system's attractor defined by a hidden Markov model, the number of time steps needed between entering and exiting is increased upon adding a random additive perturbation of moderate amplitude. Naturally, adding noise of large amplitude eventually decreases the lifetimes of the regime, so we call the effect of noise-induced enhancement of the regime lifetime for intermediate noise strength \emph{stochastic inertia}.

The first observation of stochastic inertia that we are aware of was made by Franaszek \cite{franaszek1991influence}. In numerical experiments, they considered the escape rate from a bounded set in chaotic discrete-time dynamical systems near a boundary crisis. They observed that the escape rate decreased when adding intermediate noise to the dynamics. 
%Naturally, adding large noise to the system drastically increases the escape rate since the noise term becomes dominant. 
A similar effect has been observed in \cite{altmann2010noise}, and the phenomenon has been quantified analytically in~\cite{reimann1994suppression}; see also~\cite[sections~4.1.2, 4.2.3, and 4.3.2]{lai2011transient}.
Stochastic inertia in the CdV system was first noted by \cite{kwasniok2014enhanced} along with a similar behavior in the Lorenz '63 system~\cite{lorenzDeterministicNonperiodicFlow1963}.  Importantly, the flow configuration associated with the regime considered in the CdV system vaguely resembles an atmospheric blocking and is associated with an unstable fixed point already noted in the original formulation of the model~\cite{charneyMultipleFlowEquilibria1979}.
In \cite{dorrington2023interaction}, the authors argue that stochastic inertia is caused by a flow accelerating away from the fixed point along a single unstable dimension,
% the convex flow around an unstable ``blocking'' fixed point, 
which is a property already observable in the deterministic setting. 

The aim of this article is to introduce a numerical method to identify stochastic inertia in chaotic systems. Since numerical integration of the stochastic version of the system is not always possible in practice, we introduce a method that is able to identify stochastic inertia based solely on a single (observed or simulated) trajectory of the deterministic system. Therefore, our method is purely data-driven and does not require information about the underlying dynamics. This is particularly valuable for real-world systems where only observed data is available. The idea, put simply, is to consistently emulate the effect of noise by distributing the probabilities associated with the noisy evolution on the available trajectory. Hence, we obtain a Markov chain on the data points of the trajectory.

We present the theory and derivation of the method in Sec.~\ref{sec:meth}. In Sec.~\ref{sec:1d}, we demonstrate our method on a simple one-dimensional toy system for which stochastic inertia is observed. In Sec.~\ref{sec:3d} we apply our method to a three-dimensional toy example which does not exhibit stochastic inertia. Finally, in Sec.~\ref{sec:CdV} we turn to the more complex CdV system and show that our method succeeds at identifying stochastic inertia based only on a single deterministic trajectory. We conclude by summarizing and discussing the results and point to directions of further research in Sec.~\ref{sec:fin}.

\section{Theory \& Methods}\label{sec:meth}

\subsection{Regime lifetimes}\label{sec:regime_lifetimes}
We consider dynamical systems in $\R^d$ whose trajectories transition between regimes in a seemingly unpredictable manner. The most prominent examples of such systems are those that admit a strange attractor, e.g.~the Lorenz63 system \cite{lorenzDeterministicNonperiodicFlow1963}. 

Consider an ODE
\begin{equation}\label{eq:ODE}
    \Dot{x} = v(x),
\end{equation}
for a (sufficiently) smooth vector field $v$. 
Let $\Phi^t: \R^d \to \R^d$ be the solution operator of the ODE \eqref{eq:ODE} given by
\begin{equation}
    \Phi^0 x = x \qquad \text{and} \qquad \partial_t \Phi^t x = v(\Phi^t x).
\end{equation}
An attractor of the system is a closed measurable set $A\subset\R^d$ that is invariant under $\Phi$, i.e.~$(\Phi^t)^{-1} A = A$, such that for almost every initial condition $x_0$, the distance of $\Phi^t x_0$ to the set $A$ converges to~$0$.
The mathematical theory of dynamical attractors is extremely deep, and we do not aim to provide a full technical account. For a rigorous introduction to dynamical attractors and ergodic theory, we refer, e.g., to~\cite{katok1995introduction, young2002srb}. 

Associated with the set $A$ is a measure $\mu_\Phi$ whose support comprises the attractor, i.e., $\text{supp}(\mu_\Phi) = A$. This measure is invariant under the dynamics, that is
\begin{equation}
    \mu_{\Phi}(E) = \mu_{\Phi}((\Phi^t)^{-1} E),
\end{equation}
for all $t\geq 0$ and all Borel-measurable sets $E \subset \R^d$. The characteristic property of a measure $\mu_\Phi$ representing an attractor is that the long-term statistical behavior of almost all initial conditions with respect to Lebesgue measure is governed by the measure~$\mu_{\Phi}$. % We call such a measure $\mu$ a \emph{physical measure}. 
More precisely, for any observable $f:\R^d \to \R$ and Lebesgue-almost all initial conditions $x_0\in \R^d$, we find
\begin{equation}\label{eq:Birkhoff}
    \lim_{t\to \infty} \frac{1}{t} \int_0^t f(\Phi^s x_0) ds = \int_{\R^d} f(x) \mu_{\Phi}(dx). 
\end{equation}
Invariant measures with this property are called \emph{physical measures}. A refined concept of physical measures supported on attractors are so-called \emph{SRB measures}. We refer to \cite{young2002srb} for further details. For attractors such as the Lorenz attractor, or the attractor of the Charney--deVore system (see Sec.~\ref{sec:CdV}), the set $A$ is a compact set with a fractal dimension that is smaller than the dimension of the state space. Thus, the attractor captures the asymptotic, low-dimensional dynamics of a potentially higher-dimensional system.

Let $M\subset \R^d$ be an open, bounded set (called ``regime'' from now on) of the state space that intersect the attractor, in the sense that~$\mu_\Phi(M) >0$. By \eqref{eq:Birkhoff} applied to the indicator function $\1_M$, the trajectory of almost every initial condition visits the regime $M$ infinitely often. The time between an entrance to $M$ and the subsequent escape from $M$ is called \emph{lifetime} of $M$. Let $x_t$ be a trajectory of the ODE \eqref{eq:ODE}. Without loss of generality, we assume $x_0 \notin M$, and recursively define the entrance and escape times of the regime~$M$. Set $\tau_0^\text{out} = 0$ and define for all $k\geq 1$
\begin{nalign}\label{eq:tau_in_out}
    \tau^\text{in}_k &= \inf\{ t>\tau_{k-1}^\text{out} \mid x_t \in M \},\\
    \tau_k^\text{out} &= \inf\{ t> \tau^\text{in}_k \mid x_t \notin M \}.
\end{nalign}
The differences $\tau_k^\text{out} - \tau_k^\text{in} \geq 0$ are the lifetimes of the regime $M$. In general, the trajectory $x_t$ enters and exits $M$ in an aperiodic manner, such that the lifetimes of $M$ can take a range of values. Since the attractor is ergodic, the lifetimes of $M$ have well-defined long-term statistics that are the same for almost every initial condition with respect to Lebesgue measure. In particular, there is an average regime lifetime, given by the limit
\begin{equation}\label{eq:det_regime_lifetime}
    E_M \coloneqq \lim_{K \to \infty} \frac{1}{K}\sum_{k=1}^K (\tau_k^\text{out} - \tau_k^\text{in}).
\end{equation}

We intend to study how the distribution of the lifetimes of $M$ changes when adding small additive noise to the ODE \eqref{eq:ODE}; in other words, we analyze an SDE of the form
\begin{equation}\label{eq:SDE}
    dx_t^\sigma = v(x_t^\sigma)\, dt + \sigma dW_t.
\end{equation}
As \emph{stochastic inertia} we describe the observed phenomenon that the addition of small noise increases the average regime lifetimes. As we will see, for an SDE of this form, rigorously defining average regime lifetimes becomes problematic. Therefore, also a rigorous definition of stochastic inertia is non-trivial. 

For $x_0^\sigma \notin M$, we can define the time of the first entrance to $M$
\begin{equation}
    \tau^\text{in} = \inf \{t>0 \mid x_t^\sigma \in M\},
\end{equation}
which is now a random variable.
Likewise, for $x_0^\sigma \in M$, we can define the random time of the first exit from $M$,
\begin{equation}
    \tau^\text{out} = \inf\{t>0 \mid x_t^\sigma \notin M\}.
\end{equation}
However, a recursive definition of $\tau_k^\text{in}$ and $\tau_k^\text{out}$, like in \eqref{eq:tau_in_out}, is problematic due to the nature of Brownian motion. For any time $t$ at which $x_{t}^\sigma$ is on the boundary of $M$, the trajectory enters and exits $M$ an infinite number of times in any time interval around $t$. Hence, the recursion would get stuck in the sense that $\tau_k^\text{in} = \tau_{k+1}^\text{out} = \tau_{k+1}^\text{in}$ and so on. Consequently, the average regime lifetime, as defined in \eqref{eq:det_regime_lifetime}, would equal~$0$. 

To capture the effect of prolonged regime lifetimes in an SDE, we do not work with the recursive definition \eqref{eq:tau_in_out}, and instead study the dependence of $\tau^\text{out}$ on the noise strength. To highlight this dependence in the notation, for $x_0^\sigma = x \in M$, define 
\begin{equation}\label{eq:def_tau}
    \tau^\sigma(x) \coloneqq \inf\{t>0 \mid x_t^\sigma \notin M\}.
\end{equation}
This is known as the \emph{escape time}, also called first-passage time. We emphasize that $\tau^\sigma(x)$ is a random variable, whose distribution and mean depends on the noise parameter~$\sigma$.

We note that continuous-time white noise is an abstracted model of random motion and the ill-posed nature of the average regime lifetimes can be seen as a technical artifact that plays little to no role in real-world systems. Similarly, any discrete numerical approximation of a trajectory of the SDE~\eqref{eq:SDE} naturally  enters and exits $M$ only a finite number of times. Hence, an increase in the regime lifetimes due to additive noise can be numerically observed, although an analytic quantification has yet to be developed. In the following sections, through time-stepping schemes, we treat all occurring systems as time-discrete, thus avoiding this obstacle.

By virtue of numerical experiments, we demonstrate that stochastic inertia in $M$ is linked to a noise-induced increase in the average escape times~$\E[\tau^\sigma(x)]$ for a large portion of states $x\in M$. We propose a purely data-driven method of identifying and analyzing the escape times~$\tau^\sigma$, and compare them to regime lifetimes. We emphasize that for our method no knowledge of the underlying dynamics, i.e.~of the velocity field $v$, is needed. The ``true'' dynamics is only used to verify our method. 

\subsection{Markov chain model}
\label{ssec:Markov_model}

Given a sufficiently long trajectory of a dynamical system, we construct a Markov chain that simulates how the dynamics would behave if we added small additive noise. 

Let $\{X_i \mid i=1, \hdots, N+1\}$ be a long trajectory of a dynamical system with time steps of size~$\Delta t$. It is important that the time between the samples $X_i$ and $X_{i+1}$ is the same for every~$i$. In view of the ODE \eqref{eq:ODE}, this means that~$X_{i+1} = \Phi^{\Delta t}X_i \approx X_i + v(X_i) \Delta t$. Define the set of the first $N$ points, $\mathcal{X} = \{X_i \mid i=1, \hdots, N\}$. 
Let $\mathcal{M}:= M \cap \mathcal{X} \subset \mathcal{X}$ be the set of those trajectory points that lie in the regime $M$ of interest. The set $\mathcal{M}$ is described by an index set $I\subset \{1, \hdots, N\}$ such that~$\mathcal{M} = \{X_i \mid i\in I\}$. 
We will construct a Markov chain with state space~$\mathcal{X}$ that is an approximation of the SDE~\eqref{eq:SDE} with noise strength~$\sigma$ in the weak sense; cf.~\eqref{eq:Koopman_equiv_first} below. Then, we compute the expected escape times with respect to this Markov chain for each state~$x\in \mathcal{M}$. We denote a trajectory of the Markov chain by~$\{y_n^\sigma\}_{n\in \N}$. The transition probabilities are collected in a matrix $P^\sigma\in \R^{N\times N}$ such that
\begin{equation*}
    P_{ij}^\sigma = \P(y_{n+1}^\sigma = X_j \mid y_n^\sigma = X_i).
\end{equation*}
We emphasize that the $\sigma$ in $P^\sigma$ is a superscript and not an exponent. When taking $P^\sigma$ to a power $k$, we write~$(P^\sigma)^k$.
The basis for defining the transition probabilities $P_{ij}^\sigma$ of our Markov chain is the Euler--Maruyama scheme \cite[section~9.1]{kloeden1999} for the SDE \eqref{eq:SDE}; cf.~\cite[Sec.~3.2]{koltai2018large} for a similar construction. The scheme is given by 
\begin{equation}\label{eq:Euler_Maruyama}
    z_{n+1}^\sigma = z_n^\sigma + v(z_n^\sigma) \Delta t  +  \sigma \sqrt{\Delta t} \, \xi_n,
\end{equation}
where $\{\xi_n\}_{n\in\N}$ is an iid sequence of standard-normally distributed $d$-dimensional random variables and $z_{n+1}^{\sigma}$ is the numerical approximation of a realization of the system at time $\Delta t$ when starting at time 0 in~$z_n^{\sigma}$. Our goal is to construct $P_{ij}^\sigma$ such that the Markov chain $\{y_n^\sigma\}_{n\in \N}$ in some sense approximates the Markov chain $\{z_n^\sigma\}_{n\in \N}$. The key difficulty for such an approximation is that $y_n^\sigma \in \mathcal{X}$ lives on a finite state space, while $z_n^\sigma\in \R^d$ can take any value in $\R^d$. Hence, a direct approximation in terms of the transition probabilities is impossible. Instead, we aim to approximate the dynamics $z^\sigma$ on the level of \emph{transfer operators}. 

The \emph{Koopman operator} $U$ of the Euler--Maruyama scheme \eqref{eq:Euler_Maruyama} is the linear operator on the space of bounded observables $f:\R^d \to \R$ given by
\begin{equation}
    [U f] (z) = \E[f(z_{n+1}^\sigma ) \mid z_n = z].
\end{equation}
The Koopman operator applied to an observable gives the expected value of the observable after one time step. The adjoint of the Koopman operator is called the \emph{Perron--Frobenius operator} and maps distributions of particles forward in time. For a discussion of transfer operators, we refer to~\cite{lasota2013chaos}. 
% Both the Koopman and the Perron--Frobenius operator contain all the dynamical and stochastic information of a stochastic system.

The Koopman operator of the Markov chain $y^\sigma$ on the finite state space $\mathcal{X}$ is simply given by the transition probability matrix $P^\sigma$. Indeed, let $\mathbf{f}\in \R^N$ be a column vector on $\mathcal{X}$, interpreted as an observable~$\mathbf{f}:\mathcal{X} \to \R$. Then, 
\begin{equation*}
    [P^\sigma \mathbf{f}]_i = \E[\mathbf{f}(y_{n+1}^\sigma) \mid y_n^\sigma = X_i],
\end{equation*}
is the expected value of $\mathbf{f}$ after one step of the Markov chain.
Similarly, let $\mu \in \R^N$ be a row-vector that is a probability distribution on~$\mathcal{X}$. Then $\mu P^\sigma$ is the distribution after one step of the Markov chain $y^\sigma$. This corresponds to the Perron--Frobenius operator. In short, distributions are multiplied with $P^{\sigma}$ from the left, and observables are multiplied with $P^{\sigma}$ from the right.

Let $f:\R^d \to \R$ be a bounded, measurable observable on $\R^d$. By restriction, this defines an observable $\mathbf{f}$ on $\mathcal{X}$ given by $\mathbf{f}_i=f(X_i)$. Our goal is to construct the transition probability matrix $P^\sigma$ such that $y^\sigma$ approximates $z^\sigma$  on the level of transfer operators in the sense that
\begin{equation}\label{eq:Koopman_equiv_first}
    [P^\sigma \mathbf{f}]_i \approx [U f](X_i).
\end{equation}

The transition probabilities of the Euler--Maruyama scheme \eqref{eq:Euler_Maruyama} are governed by the Gaussian kernel of a normally distributed~$\xi$. We compute
\begin{align}
    [U f](z) &= \int_{\R^d} f\big(z + v(z) \Delta t + \sigma \sqrt{\Delta t} \xi\big)  \, (2\pi)^{-\frac{d}{2}} \exp\left(-\frac{\norm{\xi}^2}{2}\right) \, \mathrm{d} \xi.
\end{align}
For $z=X_i$, we have the approximation $z + v(z) \Delta t \approx X_{i+1}$. Additionally, we can substitute $\sigma \sqrt{\Delta t}\xi = \zeta$, to arrive at
\begin{equation*}
    [U f](X_i) \approx\int_{\R^d} f\big( X_{i+1} + \zeta\big)\, (2\pi)^{-\frac{d}{2}}\exp\left(-\frac{\norm{\zeta}}{2\sigma^2 \Delta t}\right) \, \mathrm{d} \zeta.
\end{equation*}
Our goal is to construct the transition probability matrix $P^\sigma$ in such a way that 
\begin{equation}\label{eq:Koopman_equiv}
    \sum_{j=1}^N P_{ij}^\sigma \, f(X_j)  \approx \int_{\R^d} f\big( X_{i+1} + \zeta\big)\, (2\pi)^{-\frac{d}{2}}\exp\left(-\frac{\norm{\zeta}}{2\sigma^2 \Delta t}\right) \, \mathrm{d} \zeta.
\end{equation}

We construct $P^\sigma$ as the concatenation of two transition probability matrices $P^\sigma = T D^\sigma$, where
\begin{equation}
    P^\sigma \in \R^{N\times N}, \qquad T\in \R^{N \times (N+1)}, \qquad D^\sigma \in \R^{(N+1)\times N}. 
\end{equation}
Since $P_{ij}^\sigma$ is the probability of going from $X_i$ to $X_j$, the concatenation $TD^\sigma$ corresponds to first applying $T$ and then $D^\sigma$. The transition probability matrix $T$ models the deterministic step $X_i \to X_{i+1}$. Since the point $X_N$ is mapped to $X_{N+1}$, the matrix is of size $N\times (N+1)$. The matrix $T$ is given by
\begin{equation}
    T_{ij} = \begin{cases}
        1, \quad&j=i+1, \\
        0, \quad&j \neq i+1.
    \end{cases}
\end{equation}
The transition probability matrix $D^\sigma$ then models diffusion by strength $\sigma$. Since any point, including $X_{N+1}$, has to be mapped to $\mathcal{X}$, the size of this matrix is $(N+1)\times N$. In view of \eqref{eq:Koopman_equiv}, the matrix $D^\sigma$ has to be constructed in such a way that
\begin{equation}\label{eq:Koopman_equiv_final}
    \sum_{j=1}^N P_{ij}^\sigma \, f(X_j)  = \sum_{j=1}^N f(X_j) \, D^\sigma_{i+1, j} \approx \int_{\R^d} f\big( X_{i+1} + \zeta\big)\, (2\pi)^{-\frac{d}{2}}\exp\left(-\frac{\norm{\zeta}}{2\sigma^2 \Delta t}\right) \, \mathrm{d} \zeta.
\end{equation}

The construction of such a matrix $D^\sigma$ is exactly given by the diffusion map method~\cite{coifman2006diffusion}. This approach has been employed in several cases of approximating linear operators related to the dynamics~\cite{berry2015nonparametric,froyland2021spectral,froyland2024revealing}. The construction of the diffusion matrix $D^\sigma$ closely follows~\cite{coifman2006diffusion}. We write out the definition of $D^\sigma$ but refer to the above mentioned paper(s) for a deeper derivation of the formulae. 

Define the similarity matrix
\begin{equation} \label{eq:simiMat}
    K^\sigma_{ij} \coloneqq \exp\left(-\frac{\norm{X_i - X_j}^2}{2\sigma^2 \Delta t} \right), \qquad i \in \{1, \hdots, N+1\}, \quad j\in \{1, \hdots N\}.
\end{equation}
We pre-normalize $K^\sigma_{ij}$ to account for differences in the density of the point cloud $\mathcal{X}$,
\begin{equation}
    K_i^\sigma \coloneqq \sum_{j=1}^N K^\sigma_{ij}, \qquad \hat{K}^\sigma_{ij} \coloneqq \frac{K^\sigma_{ij}}{K_j^\sigma}.
\end{equation}
The diffusion matrix $D^\sigma$ is obtained by row-normalization such that it becomes a row-stochastic matrix:
\begin{equation*}
    \hat{K}_i^\sigma \coloneqq \sum_{j
    =1}^N \hat{K}^\sigma_{ij}, \qquad D_{ij}^\sigma \coloneqq \frac{\hat{K}^\sigma_{ij}}{\hat{K}_i^\sigma}.
\end{equation*}
The matrix $D^\sigma$ is a transition probability matrix that models the addition of $\sigma \sqrt{\Delta t} \, \xi_n$ to a point in $\mathcal{X}$. It is constructed such that \eqref{eq:Koopman_equiv_final} is a faithful approximation under the assumption that the point cloud $\mathcal{X}$ samples the space $\R^d$ sufficiently densely. Note that \eqref{eq:Koopman_equiv_final}, and thereby \eqref{eq:Koopman_equiv_first}, cannot possibly be good approximations for all observables $f:\R^d \to \R$. In particular, if $f$ is non-zero but $f(X_i) = 0$ for all $X_i\in \mathcal{X}$, then $P^\sigma \mathbf{f} = 0$ while $U f$ is not close to 0. For a more detailed discussion of the necessary assumptions and limitations of this method, we refer to Section \ref{sec:lim}.

Computing all pairwise distances in \eqref{eq:simiMat} is computationally expensive for large $N$ and leads to a dense matrix of size $(N+1)\times N$. Note that for pairs $X_i, X_j \in \mathcal{X}$ for which $\norm{X_i - X_j}^2$ is of an order larger than $2\sigma^2 \Delta t$, the value $K_{ij}^\sigma$ is extremely small and can be neglected. Therefore, we calculate the distances needed for \eqref{eq:simiMat} with a nearest neighbour search using a \texttt{kd-tree} as implemented in \texttt{scipy} \cite{2020SciPy-NMeth} with a cut-off radius of $3\sqrt{2\sigma^2 \Delta t}$. This leads to sparse matrices $K^\sigma$ and $D^\sigma$ which significantly improves the computational efficiency and decreases storage requirements for large~$N$. %An exception is of course the one-dimensional case, where distances are trivially available. \rc{I asked Henry about the 1D case.}
Note that if $X_N$ is the only point of $\mathcal{X}$ within the cut-off radius of the point $X_{N+1}$, the resulting Markov chain will feature a closed single-node loop. To avoid this, we extend the cut-off radius just enough to include a second point in those cases. 

\begin{remark}
\label{rem:noisy_basetraj}
    We note that this method can analogously be used when the $X_i$ are samples from a trajectory with fixed \emph{nonzero} noise, say~$\sigma_0$. The matrix $P^{\sigma}$ in this case approximates the transfer operator of the system with noise~$\sqrt{\sigma_0^2 + \sigma^2}$. Since natural systems are believed to be intrinsically noisy, the method thus allows to analyze the effect of \emph{increasing} noise in the system.
\end{remark}

\begin{remark}
    Another approach approximating the dynamics of a continuous-state system by a Markov chain is \emph{Ulam's method}~\cite{Ulam60}, which defines transition probabilities between elements of the state space's partition. 
    %The method we consider has the advantage over 
    Note that, in contrast to Ulam's method, our approach does not require the choice of a partition, making, in particular, the inclusion of variable noise strengths $\sigma$ computationally more straightforward.
\end{remark}

\subsection{Regime lifetimes of the Markov chain}

Having constructed the Markov operator $P^\sigma$ for each $\sigma>0$, we compute the average escape times from $\mathcal{M}$ under $P^\sigma$. Here and in the following we assume that the Markov chain is irreducible, meaning that there is no true subset $J \subsetneq \{1,\ldots,N\}$ of indices which is invariant with probability one.  
This assumption is justified in the setting we are interested in due to the mixing properties of chaotic attractors.
Recall that $I \subset \{1,\hdots, N\}$ is the set of indices of points $X_i$ that belong to $\mathcal{M}$. Let $Q^\sigma$ be the submatrix of $P^\sigma$ with indices $I\times I$. Hence, $Q^\sigma$ is of size $|I|\times |I|$. To avoid unnecessary complications, we index $Q^\sigma$ with the indices of $I$ instead of by $1, \hdots, |I|$. We also write $Q^\sigma \in \R^{I\times I}$. Later on, we use the notation $\R^I$ to denote the space of $|I|$-dimensional vectors that are indexed by the set $I$. 

The matrix $Q^\sigma$ models the Markov chain inside $\mathcal{M}$ and is substochastic because there is a positive probability of leaving $\mathcal{M}$ by virtue of the Markov chain's irreducibility and $I \neq \{1,\ldots,N\}$. This also means that $Q^\sigma$ has spectral radius strictly less than~1. Assuming that $y_0^\sigma = X_i \in \mathcal{M}$, the probability of remaining inside of $\mathcal{M}$ for $k$ steps under the Markov chain $P^\sigma$ is given by
\begin{equation}
    \sum_{j\in I}\left[(Q^\sigma)^k\right]_{ij} = \left[(Q^\sigma)^k\1 \right]_i,
\end{equation}
where $\1$ is a vector with each entry equal to~$1$. Hence, the average number of steps until leaving $\mathcal{M}$, when started in $X_i$, is given by~\cite[Theorem 1.3.5]{norris1998markov},
\begin{align*}
    \E[\text{no.\ of steps to leave $\mathcal{M}$} &\mid \text{start at }X_i] \\ 
    &= \sum_{k=1}^{\infty} k\, \P[\text{exactly $k$ steps to leave $\mathcal{M}$} \mid \text{start at }X_i] \\
    \text{(by counting terms)} &= \sum_{k=0}^{\infty} \sum_{l=k+1}^{\infty} \P[\text{exactly $l$ steps to leave $\mathcal{M}$} \mid \text{start at }X_i] \\
    &= \sum_{k=0}^{\infty} \P[\text{stay $k$ steps before leaving $\mathcal{M}$} \mid \text{start at }X_i] \\
    &= \sum_{k=0}^\infty \left[(Q^\sigma)^k\1 \right]_i = \left[ (\Id- Q^\sigma)^{-1} \1 \right]_i.
\end{align*}
The formula for the power series
% as the resolvent of $Q^\sigma$
is valid since $Q^\sigma$ is substochastic and therefore has spectral radius strictly less than~$1$. Define the vector $\theta^\sigma \in \R^{I}$ by
\begin{equation}\label{eq:def_theta}
    \theta^\sigma \coloneqq \Delta t (\Id - Q^\sigma)^{-1}\1.
\end{equation}
This vector approximates the escape times \eqref{eq:def_tau} of $M$
such that for each $i\in I$,
\begin{equation}
    \theta^\sigma_i \approx \E[\tau^\sigma(X_i)].
\end{equation}

We solve the linear system $(\Id - Q^{\sigma})x=\1$ to compute $\theta^{\sigma}$ numerically. For dense matrices, a direct Linear Algebra PACKage (LAPACK) solver is used, while for sparse matrices we first attempt a sparse LU factorization. If this is unsuccessful, we use an iterative Generalized Minimal Residual (GMRES) method with incomplete LU (ILU) preconditioning if necessary, falling back to a dense solve only if other methods fail to adapt to varying sparsities of $Q^{\sigma}$. We rely on the implementation by \texttt{scipy} \cite{2020SciPy-NMeth}.

Using the escape times $\theta^\sigma$, we can compute the average regime lifetimes of $\mathcal{M}$ under the Markov chain $P^\sigma$. To compute the average lifetime, we need two types of information. How likely is it to enter $\mathcal{M}$ at a point $X_i \in \mathcal{M}$, and what is the expected escape time from it when starting in~$X_i$? The latter is given by the vector $\theta^\sigma$. The entry probabilities are given by
\begin{equation}\label{eq:conditional_probabilities}
    \nu_i^\sigma \coloneqq \P(y_{n+1}^\sigma = X_i \mid y_n^\sigma \notin \mathcal{M} \: \land y_{n+1}^\sigma \in \mathcal{M} ).
\end{equation}
These probabilities are computed in the following way: Start with the stationary distribution outside of $\mathcal{M}$, apply $P^\sigma$ once, restrict the result to $\mathcal{M}$ and normalize. Let $\mu^\sigma \in \R^N$ be the stationary distribution of $P^\sigma$, i.e.~the (by irreducibility unique) vector such that $\mu^\sigma P^\sigma = \mu^\sigma$. The (not normalized) stationary distribution outside of $\mathcal{M}$ is denoted by $\mu^\sigma_{\mathcal{M}^c}$ and given by
\begin{equation}
    \left(\mu^\sigma_{\mathcal{M}^c}\right)_i \coloneqq \begin{cases}
        \mu^\sigma_i, \qquad & i \in I, \\
        0, \qquad & i \notin I.
    \end{cases}
\end{equation}
The formula for the conditional probabilities \eqref{eq:conditional_probabilities} reads 
\begin{equation}\label{eq:def_nu}
    V_i^\sigma \coloneqq [ \mu^\sigma_{\mathcal{M}^c} P^\sigma]_i, \qquad \nu_i^\sigma = \frac{V_i^\sigma}{\sum_{j\in I} V_j^\sigma}.
\end{equation}
We note that it is not necessary to normalize $\mu^\sigma_{\mathcal{M}^c}$, since $\nu^\sigma$ is normalized in the end anyway. Finally, the expected lifetimes of the regime $\mathcal{M}$ are given by
\begin{equation} \label{eqn:expRegLif}
    E^\sigma \coloneqq \langle\nu^\sigma, \theta^\sigma\rangle,
\end{equation}
where we take the scalar product over two vectors in $\R^I$.

\subsection{Properties of the method} \label{sec:lim}
One of the key advantages of the Markov chain method is that it is purely data-driven and, unlike the Monte Carlo approximation of average regime lifetimes, does not require information on the underlying dynamics or, in an observational setting, does not require data from noisy cases.
Nevertheless, its accuracy is also limited by attributes of the available trajectory data. In the following, we are making the reader attentive to some other important aspects that will partly be relevant in the numerical explorations in Sec.~\ref{sec:app}.

\paragraph{Resolution.}

\textbf{a)} The Markov chain method is \emph{global} in the sense that it utilizes interactions between states due to virtual noise. For its accuracy, a sufficient global resolution by the points $X_i$ is necessary. For the expected escape times $\theta^\sigma$, the regime should be well resolved; for the entry probabilities $\nu^\sigma$, even the whole attractor. For a low-dimensional attractor, this seems computationally possible. However, if the points $X_i$ stem from a trajectory, they sample the attractor with respect to the invariant measure~$\mu_{\Phi}$. This measure can be arbitrarily nonuniform, resulting in extremely low sample densities in certain regions. It is unclear whether an inferior resolution of these ``dynamically less relevant'' regions affects the performance of the Markov chain method.

\textbf{b)} Literature related to diffusion maps~\cite{coifman2006diffusion,berry2015nonparametric} shows that the matrix $D^{\sigma}$ from Sec.~\ref{ssec:Markov_model} models diffusion on the attractor well only within a certain range of $\sigma$, which depends on the local point density and curvature of the attractor. A too small $\sigma$ perceives a ``disconnected'' manifold, while a too large $\sigma$ ignores its geometric features. The ``good'' region of $\sigma$ values is found as a nonconstant linear region in the log-log representation of $\sum_{i,j}K^\sigma_{ij}$ versus~ $\sigma$. The slope of this curve in the good region approximates the dimension of the manifold (cf.~also~\cite{schoellerAssessingLagrangianCoherence2025} for an application). In our examples below, this is utilized in Figures~\ref{fig:sigmaloglogtoy3d} and~\ref{fig:sigmaloglogcdv}.

\textbf{c)} Finally, the time span between consecutive samples $X_i$ and $X_{i+1}$ should not be too large (depending on the nonlinearity of the deterministic flow); otherwise the Euler--Maruyama approximation in Eqn.~\eqref{eq:Euler_Maruyama} deteriorates. Since the Markov chain model is constructed to approximate the Euler--Maruyama scheme, large time steps in turn constraints the accuracy of our model.

\textbf{d)} Increasing the number of steps has the benefit of potentially resolving the attractor better, but comes at the cost of an increased number of states.
While the computational effort of computing the expected escape times $\theta^\sigma$ scales with the size of the regime $\mathcal{M}$ (a linear system has to be solved), the computational effort of computing the entry distribution $\nu^\sigma$ scales with the size of the entire trajectory $\mathcal{X}$ (an eigenvalue problem has to be solved), and computing $\nu^\sigma$ quickly becomes unfeasible without additional approximation methods. 

\paragraph{Restriction to deterministic attractor.}

\textbf{e)} Perturbing the deterministic dynamics by isotropic white noise will perpetually ``kick'' the trajectory off the attractor. One expects that, due to attraction, the noisy trajectory will, with overwhelming probability, stay within a distance to the attractor that is commensurate with the noise strength. The formalization of this intuition is beyond the scope of our work; we merely refer to a possible approach~\cite{blank1986stochastic}. Nevertheless, this change in the ``effective attractor'' may necessitate redefining the regime (due to, for instance, some physical considerations). The Markov chain method cannot account for this change, leaving an intrinsic bias in its results.

\textbf{f)} Even if the effective attractor changes only negligibly, the interactions between the flow towards the attractor and perturbations transversal to the attractor can have a relevant effect. This effect cannot be modelled entirely via dynamical information \emph{on} the attractor. However, the Markov chain method is restricted to noise that acts tangentially\footnote{Here, we assume that there is a well-defined tangential space of the attractor. If not, for instance because the attractor is not a smooth manifold, the space spanned by the neighborhood of any point of the attractor might still be strictly smaller than the entire~$\R^d$.} to the attractor and cannot incorporate noise that drives the trajectory away from the attractor. 
This discrepancy arises only if the attractor has a lower dimension than the full state space, but in this case, the effect may be severe.

\section{Application \& Results}\label{sec:app}

\subsection{Stochastic Inertia in a one-dimensional toy system}\label{sec:1d}

We first demonstrate the methodology developed in Sec.~\ref{sec:meth} by applying it to a one-dimensional toy model. Naturally, it is not possible to have an attractor in one dimension that is not simply a fixed point. Instead, imagine an attractor in higher dimensions that contains a fixed point with one unstable direction while all other directions are stable. Trajectories approach the fixed point through one of the stable directions, spend some time near the fixed point, and ultimately leave through the unstable direction. After spending some time in the chaotic region of the attractor, the trajectory returns in a seemingly random manner to the fixed point. Thus, in the vicinity of this fixed point, the dynamics is approximated by the one-dimensional dynamics in the unstable direction with random reinsertions. For example, such a setting arises in the blocking regime in the Charney--deVore model; see Sec.~\ref{sec:CdV}. Heuristically, this blocking regime has previously been treated as a one-dimensional system in \cite{dorrington2023interaction} to explain the occurrence of stochastic inertia.  

\begin{figure}[]
    \centering
    \includegraphics[width=.75\textwidth]{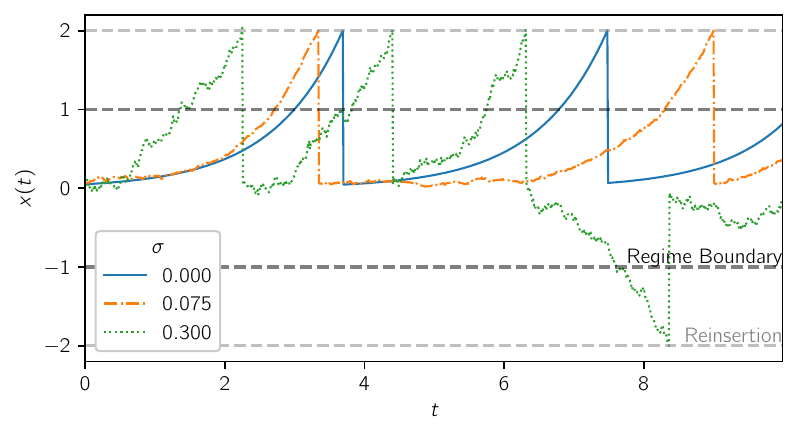}
    \caption{Example trajectories initialized in $x_0 = 5\times 10^{-2}$ for the one-dimensional toy model \eqref{eq:1d_toy} for three different noise strengths $\sigma$, $10^{3}$ steps of $\Delta t=10^{-2}$ and a reinsertion rule with $\epsilon = 5\times 10^{-2}$. Also shown are the boundaries defining the regime $M$ and the system's boundaries at whose crossing reinsertion takes place.}
    \label{fig:Toy}
\end{figure}

\begin{figure}[]
    \centering
    \includegraphics[width=.75\textwidth]{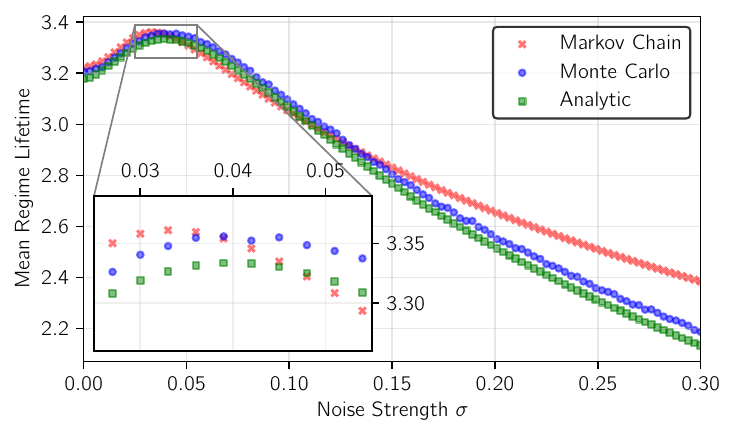}
    \caption{Mean regime lifetimes for each noise strength $\sigma$
obtained from the Monte Carlo method, the Markov chain approach or the analytic solution. The former two methods both use time steps of $\Delta t = 5\times 10^{-2}$ starting from $X_1 = z_0 = 10^{-4}$ for comparability, but the Monte Carlo method uses $5 \times 10^6$ steps containing $\approx 10^4$ reinsertions in contrast to $2 \times 10^4$ steps with $250$ reinsertions for the Markov Chain method. Reinsertion is done according to Eqn.~\eqref{eq:reinsertion} with $\epsilon=5\times 10^{-2}$. The analytic curve is calculated by averaging solutions of Eqn.~\eqref{eq:analytic_1d} over $250$ values of $R_k$. 
% The distribution of reinsertions is provided as Fig.~\ref{fig:Toy1dReinsert} in the appendix.
}
    \label{fig:TOLifetInc}
\end{figure}

We consider a toy model of an unstable fixed point in one dimension. Once trajectories are sufficiently far away from the fixed point, we randomly reinsert them close to the fixed point. The system is given by
\begin{equation}\label{eq:1d_toy}
    \mathrm{d}x_t = x_t\,\mathrm{d}t + \sigma\,\mathrm{d}W_t,\qquad x_t\in(-2,2), \qquad M = (-1,1),
\end{equation}
where $M$ is the regime we investigate. Maintaining some space between the regime boundary and the domain boundary is required since our Markov chain method relies on some non-empty ``outside'', i.e. some non-trivial complement of the regime. Reinsertion of the system occurs upon hitting the boundary with reinsertion times 
\begin{equation}
    \tau_0=0,\qquad \tau_{k+1}=\inf\{t>\tau_k \mid |x_t|\ge 2\},
\end{equation}
and reinsertion rule
\begin{equation}
    x_{\tau_{k}} = \mathrm{sign}\big(x_{\tau_k^-}\big)\,R_k, \quad k\ge 1,
\end{equation}
where $R_k$ is drawn according to 
\begin{equation}\label{eq:reinsertion}
    R_k = \sqrt{\sum_i^3 c_i^2}\, , \qquad c_i \sim \mathcal{N}(0, \epsilon),
\end{equation}
where $\epsilon > 0$ is small. Here, $x_{\tau_k^-}$ denotes the boundary point from which the reinsertion takes place, i.e.,~$\lim_{t \nearrow \tau_k} x_t$. Hence, $R_k$ is the Euclidean norm of a small three-dimensional vector with each component drawn from a normal distribution. For the results shown in Fig.~\ref{fig:TOLifetInc}, we used $\epsilon = 5\times10^{-2}$. The reinsertion law described by Eqn.~\eqref{eq:reinsertion} models randomly approaching a fixed point in a higher-dimensional space through stable directions. In Fig.~\ref{fig:Toy}, the trajectories of Eqn.~\eqref{eq:1d_toy} for three values of $\sigma$ are shown. Though not statistically meaningful, the figure shows how the regime lifetime can be enhanced by adding intermediate noise. While the trajectory with large noise ($\sigma = 0.3$) quickly leaves the vicinity of the fixed point, the trajectory with intermediate noise ($\sigma = 0.075$) spends a long time close to the fixed point, thus producing long regime lifetimes.

For this one-dimensional system, the lifetimes of the regime $M$ are given by the times between reinsertion and escape from~$M$. We apply three methods to compute the average regime lifetimes for fixed initial conditions: The Markov chain method introduced in Sec.~\ref{sec:meth}, a Monte Carlo simulation, and the analytic solution. 
% One of the key properties of the Markov chain method is that it is purely data-driven and, unlike the Monte Carlo simulation, does not require information on the underlying dynamics or, in an observational setting, does not require data from the noisy cases.

To apply the Markov chain method introduced in Sec.~\ref{sec:meth}, we sample a long trajectory $X_1, \hdots, X_{N+1}$ of Eqn.~\eqref{eq:1d_toy} without noise, i.e.~$\sigma=0$. Note that even though $\sigma=0$, the trajectory is random since the reinsertions $R_k$ are random. In order for the method to work optimally, the sampled trajectory should explore the entire state space. To guarantee that the trajectory explores the immediate proximity of zero, we initialize at~$X_1 = 10^{-4}$. We take $N=2\times 10^4$ samples with a step size of~$\Delta t = 5\times 10^{-2}$. This results in $250$ reinsertions.
% ; a histogram of which is shown in Fig.~\ref{fig:Toy1dReinsert} in the appendix. 
Given $\sigma \geq 0$ we compute the average escape times $\theta^\sigma$, as described in~Eqn.~\eqref{eq:def_theta}. To obtain the average lifetimes of $M$ using the Markov chain method, we have to integrate $\theta^\sigma$ over the reinsertion distribution. Since the reinsertion distribution does not depend on $\sigma$, we simply take the average of $\theta^\sigma$ over all reinsertion points. Let $K\in \N$ be the number of reinsertions (in our case $K=250$) and let $X_{i_1}, \hdots, X_{i_K}$ be the reinsertion points. Then we compute the average regime lifetime by
\begin{equation}
    E^\sigma = \frac{1}{K} \sum_{k=1}^K \theta^\sigma_{i_k}.
\end{equation}

We compare the quantity $E^\sigma$ to the average regime lifetimes obtained using a Monte Carlo simulation. Given $\sigma \geq 0$, we initialize a random trajectory at $z_0^\sigma= 10^{-4}$ and compute $5\times 10^6$ steps with $\Delta t = 5\times 10^{-2}$. The numerical integration scheme is based on the following observation: Without reinsertions, the solutions of Eqn.~\eqref{eq:1d_toy} are given by
\begin{equation}\label{eq:mild_form}
    x_t^\sigma = e^t x_0^\sigma + \sigma \int_0^t e^{t-s} dW_s.
\end{equation}
Let $z_n^\sigma$ denote the numerical solution at time $n\Delta t$ and let $\xi_n \sim \mathcal{N}(0,1)$ be i.i.d. Then
\begin{equation}\label{eqn:Toy1dNum}
z_{n+1}^\sigma=
    \begin{cases}
        \mathrm{sign}(z_n^\sigma)\,R_{k}, & \text{if }|z_n^\sigma|\ge 2,\\
        e^{\Delta t} z_n^\sigma + \xi_n\,\underbrace{\sigma\,\sqrt{\frac{e^{2\Delta t}-1}{2}}}_{\eqqcolon \sigma_{\Delta t}}, & \text{if }|z_n^\sigma|<2.
    \end{cases}
\end{equation}
Scaling the Gaussian term by $\sigma_{\Delta t}$ ensures that the discrete update has the correct variance matching~Eqn.~\eqref{eq:mild_form}.  Applying an Euler--Maruyama approximation would, by contrast, be correct only asymptotically as $\Delta t \rightarrow 0$ and would produce a timestep-dependent effective noise level. Note that \eqref{eqn:Toy1dNum} allows for $|z_n^\sigma|>2$ for a single time step before reinsertion, while the SDE \eqref{eq:1d_toy} is defined on $(-2, 2)$. However, this does not diminish the faithfulness of the method. Trivially, the trajectory $X_1, \hdots, X_{N+1}$ we use for the Markov chain method is calculated with the same scheme with~$\sigma=0$. A highly parallelized, just-in-time compiled and single precision implementation of the system is provided in~\cite{schoeller_2026_19711549}.

Due to the simple nature of the SDE \eqref{eq:1d_toy}, there is an analytic expression for the average escape time for any $x\in M$. For $\sigma>0$ it is given by
\begin{equation}\label{eq:analytic_1d}
    \E[\tau^\sigma (x)] = \frac{\int_{-1}^x \exp\left(-\frac{y^2}{\sigma^2}\right) \int_{-1}^y \exp\left(-\frac{y^2 - z^2}{\sigma^2}\right) \, \mathrm{d}z\, \mathrm{d}y  }{\frac{1}{2}\sigma^2 \int_{-1}^1 \exp\left(-\frac{y^2}{\sigma^2}\right) \, \mathrm{d} y}.
\end{equation}
For a derivation of the formula, see \cite[Eqn.~(6.34)]{schuss2009theory}. Numerical evaluation of this formula is briefly described in App.~\ref{app:NumRem}. 
Computing the analytical regime lifetimes for each $\sigma$ equates to an average of $\E[\tau^\sigma (x)]$ from Eqn.~\eqref{eq:analytic_1d} over $x$ weighted by the reinsertion probability distribution. By the law of large numbers, we expect the analytical average escape times to agree with the average escape times of the Monte Carlo simulation.

The average regime lifetimes, computed using the Markov chain method, the Monte Carlo simulation, and the analytic solution are shown in Fig.~\ref{fig:TOLifetInc}. For small noise strength up to $\sigma \approx 0.1$, the three methods produce similar results, and clearly show stochastic inertia. The highest average regime lifetime is reached for $\sigma \approx 0.04$. For larger values of $\sigma$, the Markov chain method significantly overestimates the average regime lifetimes. This indicates that our method is able to identify stochastic inertia, but is not suited to compute regime lifetimes for large noise strengths (see Sec.~\ref{sec:lim}, b)). Some of the minor disagreement between the three calculation methods are perhaps attributable to the differences in the distribution of reinsertions.

\begin{figure}[]
    \centering
    \includegraphics[width=.6\textwidth]{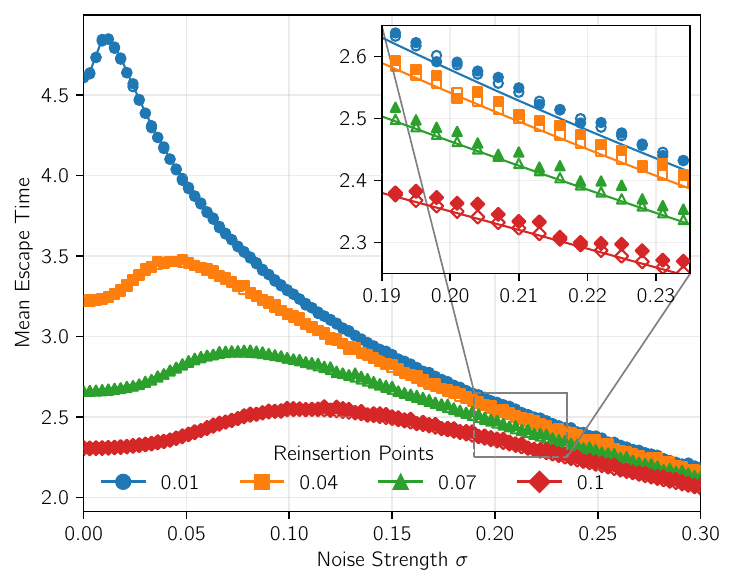}
    \caption{Mean escape times for the toy system \eqref{eq:1d_toy} initialized in various points over a range of noise strengths $\sigma$. The lines indicate analytic solutions according to Eqn.~\eqref{eq:analytic_1d}. The full markers indicate empirical mean escape times of Monte Carlo simulations of $N \approx 10^4$ realisations with time steps of $\Delta t = 10^{-2}$. The hollow markers indicate the expected mean escape times $\theta^\sigma$ obtained using the Markov chain method and a deterministic trajectory initialized in $X_1 = 10^{-4}$ and ending upon exiting the domain. This amounts to $19,807$ steps with $\Delta t = 5\times10^{-4}$. Both trajectories were produced according to Eqn.~\eqref{eqn:Toy1dNum}. The inset details the tail of the curves.}
    \label{fig:Toy_AnaVSTraj}
\end{figure}

To test the theoretical accuracy of our method, we compute the average escape times for four fixed initial conditions using the Markov chain method. Since no reinsertions are needed for this computation, we compute one long trajectory, started at $X_1 = 10^{-4}$ that ends upon leaving the domain $(-2, 2)$. To increase the density of points, we use a significantly smaller step size of $\Delta t = 5 \times 10^{-4}$ and end up with $19,807$ points. The resulting average escape times are depicted in Fig.~\ref{fig:Toy_AnaVSTraj}, and compared to the corresponding average escape times of a Monte Carlo simulation with a long trajectory of $5 \times 10^6$ steps of $\Delta t = 10^{-2}$ and the analytical solutions. It is not expected that points used as initial conditions are part of the trajectory $\mathcal{X}$ which we use for the Markov chain method. Instead, we find the closest points in $\mathcal{X}$ and show their expected escape times. They agree with the chosen initial conditions up to a relative error of $10^{-4}$ at most. The results shown in Fig.~\ref{fig:Toy_AnaVSTraj} demonstrate the behavior expected for systems featuring stochastic inertia: The average escape time is significantly higher for intermediate noise levels compared with the deterministic ($\sigma=0$) system or to high noise levels. For initial conditions closer to zero, the value of $\sigma$ that maximizes the average escape time is smaller, but the absolute increase of the escape times is of the same order for all initial conditions we tested. 

For the chosen parameters, the three methods are remarkably accurate. This indicates that the Markov chain method is able to compute expected escape times extremely well, assuming the initial trajectory is dense enough. In particular, increasing the time step $\Delta t$ makes the Markov chain method's performance for higher noise values deteriorate. This explains the incorrect tail visible in Fig.~\ref{fig:TOLifetInc}. Also, we need to keep in mind that the simulation and, in particular, the measurement of a system in as much detail as in this example is typically not feasible.

\subsection{A three-dimensional example without stochastic inertia}\label{sec:3d}
In \cite{dorrington2023interaction} a heuristic explanation of stochastic inertia around a one-dimensional unstable fixed point is given. In brief terms, their reasoning is as follows: Taking a random step towards the boundary of the regime saves less time than a random step towards the fixed point costs for reaching the regime boundary. Thereby, adding small zero-mean noise to the system increases the average escape time. Naturally, as the noise gets larger, the noise becomes dominant and trajectories leave the domain quickly, thus resulting in an increase in regime lifetimes only for intermediate noise levels. 

The reasoning from \cite{dorrington2023interaction} is specific to one-dimensional systems that accelerate towards the boundary. We consider a three-dimensional analogue of the one-dimensional system \eqref{eq:1d_toy} and show that there is no stochastic inertia in this example. We note that this is not an argument that stochastic inertia cannot occur in higher dimensional systems. Adding stable directions does not play a major role for the escape times from a regime. The CdV system considered in Sec.~\ref{sec:CdV} is six-dimensional and the blocking regime still exhibits stochastic inertia. However, the blocking regime has one dominant unstable direction and, thus, behaves locally like a one-dimensional system.

We consider a toy model of an unstable fixed point in three dimensions. Once trajectories are sufficiently far away from the fixed point, we randomly reinsert them close to the fixed point. The system is given by
\begin{equation}\label{eq:3d_toy}
    \mathrm{d}x_t = x_t\,\mathrm{d}t + \sigma\,\mathrm{d}W_t,\qquad \norm{x_t}\in(-2,2), \qquad M = \{ x\in \R^3 \mid \norm{x} < 1\},
\end{equation}
where $M$ is the regime we investigate and $\norm{\cdot}$ the Euclidean norm. Reinsertion is analogous to the one-dimensional case with reinsertion times 
\begin{equation}
    \tau_0=0,\qquad \tau_{k+1}=\inf\{t>\tau_k \mid \norm{x_t}\ge 2\},
\end{equation}
and reinsertion rule
\begin{equation}
    x_{\tau_{k}} = R_k, \quad k\ge 1,
\end{equation}
where $R_k$ is drawn according to 
\begin{equation}\label{eq:reinsertion3d}
    R_k = (c_1, c_2, c_3) , \qquad c_i \sim \mathcal{N}(0, \epsilon),
\end{equation}
with independent components and small~$\epsilon > 0$. Temporal discretisation is according to Eqn.~\eqref{eqn:Toy1dNum}. 

\begin{figure}[]
    \centering
    \includegraphics[width=.6\textwidth]{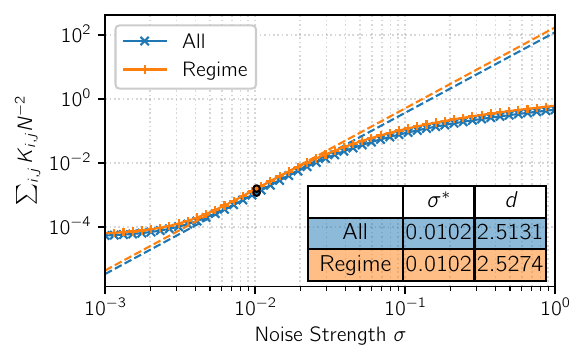}
    \caption{Normalized sum of affinities according to Eqn.~\eqref{eq:simiMat} plotted over the noise strength $\sigma$ for all points used in the analysis (``All'') as well as for points in $M$ (``Regime''). The dashed lines show the tangent at $\sigma = \sigma^{\ast}$ (indicated by open circles), the noise level maximizing the slope of the curve with the maximum slope~$d$. Both values are indicated for both sets in the table in the bottom right.}
    \label{fig:sigmaloglogtoy3d}
\end{figure}

Like in the one-dimensional example, we use the Markov chain method to compute the average regime lifetimes. The only algorithmic change that occurs when transitioning from one to multiple dimensions is that the point-wise distances needed in Eqn.~\eqref{eq:simiMat} are now calculated by a nearest neighbour search using a kd-tree, whereas in the one-dimensional case, distances were trivially available. In addition, we have to sample the three-dimensional state space sufficiently densely, or else the Markov chain method will behave just like in the one-dimensional case. To check whether the sampled trajectory $\mathcal{X}$ is sufficiently dense, we use the heuristics described in Sec.~\ref{sec:lim}, b). That is, the Markov chain method is expected to work for values of $\sigma$ in a range where the relation between $\sigma$ and $\sum_{i,j} K_{ij}^\sigma$ is linear in a log-log representation. Such a plot is shown in Fig.~\ref{fig:sigmaloglogtoy3d} for the trajectory we will use afterwards. The curve appears linear in a range of roughly $3\times 10^{-3} \leq \sigma \leq 2 \times 10^{-2}$ with the maximum slope occurring well inside our region of focus. The heuristically derived dimensionality $d_\mathcal{X}$ is around $2.5$ which falls short of the true dimension. However, samples of the same size but homogeneously distributed with respect to direction and radius exhibit dimensionalities only marginally higher (not shown here), indicating that underestimation of the dimension is a non-circumventable boundary artifact. This behavior has been observed in \cite{schoellerAssessingLagrangianCoherence2025} already. The exponentially decreasing point density by radius leads to the curve in Fig.~\ref{fig:sigmaloglogtoy3d} saturating slowly, but this is of little importance to us since we are mainly interested in the low-noise behavior. The curve for the whole trajectory and the one for only $\mathcal{M}$ are virtually the same for the same reason.

\begin{figure}[]
    \centering
    \includegraphics[width=.75\textwidth]{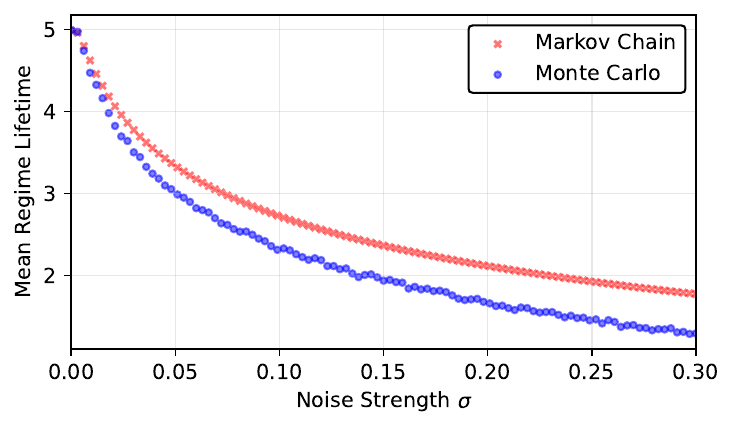}
    \caption{Mean regime life times for each noise strength $\sigma$ obtained from either the Monte Carlo method or using the Markov chain approach. Both use $2\times 10^4$ time steps of $\Delta t = 10^{-1}$ starting from $x_0 = 5\times 10^{-2} \, \1$ for comparability. The Monte Carlo trajectories feature roughly $10^2$ reinsertions with $344$ reinsertions for the Markov Chain method. Reinsertion is done according to Eqn.~\eqref{eq:reinsertion3d} with $\epsilon =5 \times  10^{-3}$. %The distribution of reinsertions for the trajectory used for the Markov Chain method is provided as Fig.~\ref{fig:Toy3dSampling} in the appendix.
    }
    \label{fig:Toy3dLifetInc}
\end{figure}

Having checked that the sampled trajectory is dense and homogeneous enough for our method to work, we compare its results with the ones obtained from Monte Carlo simulations in Fig.~\ref{fig:Toy3dLifetInc}. We use a shorter trajectory for the Markov chain method due to computational restrictions, cf.\ Sec.~\ref{sec:lim}~d). Unlike in the one-dimensional example, we do not observe stochastic inertia in the three-dimensional example. Both methods agree in this respect. As in the one-dimensional example, the Markov chain method overestimates the regime lifetimes for higher noise values -- an effect that, again, becomes more prominent when increasing~$\Delta t$. As discussed above with regard to Fig.~\ref{fig:sigmaloglogtoy3d}, this is not caused by shortcomings of the diffusion maps method, but can be rather explained by the limitations explained in Sec.~\ref{sec:lim} b).

\subsection{Stochastic Inertia in the CdV System}
\label{sec:CdV}

\paragraph{The CdV system and its attractor.}

As a final proof-of-concept, we turn to the Charney--deVore system, inspired by the work~\cite{dorrington2023interaction}. The system consists of six ordinary differential equations. They are derived by projecting a stream function field that behaves according to simplified atmospheric physics onto six Fourier modes in a Galerkin approximation. For a thorough derivation, parameter values, and analysis we refer to~\cite{deswartLoworderSpectralModels1988, crommelinMechanismAtmosphericRegime2004, dorrington2023interaction}. The system is given by 
% \rc{what are $x_1^*$ and $x_4^*$?}
\begin{align} \label{eq:CdV}
    \dot{x}_{1} &= -C(x_{1} - x_{1}^{*}) + \tilde{\gamma}_{1}x_{3}\\
\dot{x}_{2} &= -Cx_{2} + \beta_{1}x_{3} - \alpha_{1}x_{1}x_{3} - \delta_{1}x_{4}x_{6}\\
\dot{x}_{3} &= -Cx_{3} - \beta_{1}x_{2} + \alpha_{1}x_{1}x_{2} + \delta_{1}x_{4}x_{5} - \gamma_{1}x_{1}\\
\dot{x}_{4} &= -C(x_{4} - x_{4}^{*}) + \epsilon(x_{2}x_{6} - x_{3}x_{5}) + \tilde{\gamma}_{2}x_{6}\\
\dot{x}_{5} &= -Cx_{5} + \beta_{2}x_{6} - \alpha_{2}x_{1}x_{6} - \delta_{2}x_{4}x_{3}\\
\dot{x}_{6} &= -Cx_{6} - \beta_{2}x_{5} + \alpha_{2}x_{1}x_{5} + \delta_{2}x_{4}x_{2} - \gamma_{2}x_{4},
\end{align}
where $C, \tilde\gamma_i, \gamma_i, \beta_i, \alpha_i, \delta_i, \epsilon$ are model parameters and $x^*_i$ defines a background state towards which the system relaxes.

\begin{figure}[H]
    \centering
    \includegraphics[width=0.95\textwidth]{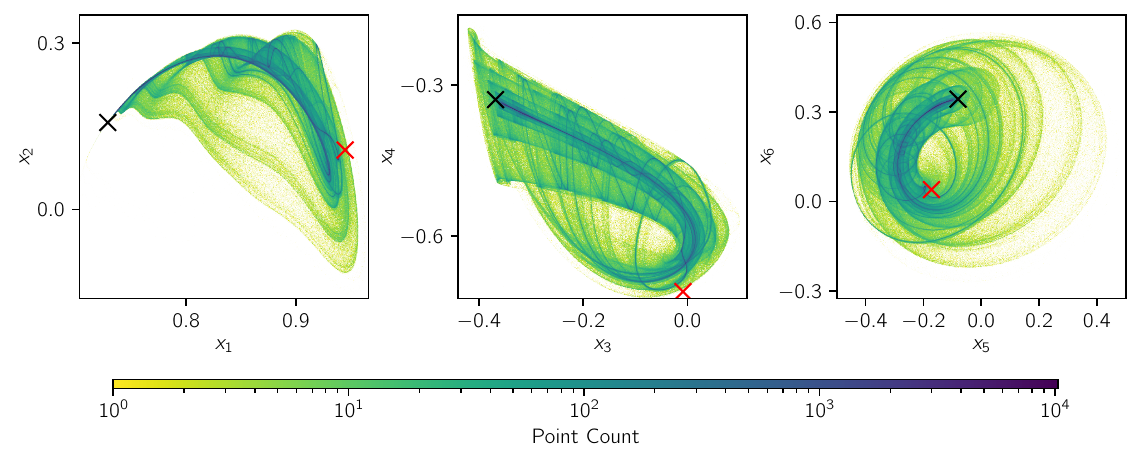}
    \caption{Two-dimensional projections of the deterministic CdV attractor based on a trajectory of $10^7$ time steps with $\Delta t=1$. The figures show two dimensional histograms with $10^3 \times 10^3$ bins uniformly partitioning the rectangular bounding box of the data.  A few outliers have been clipped. Also indicated are the ``blocking'' fixed point (black) and the ``zonal'' fixed point (red). Adapted from \cite{dorrington2023interaction}.}
    \label{fig:Hists}
\end{figure}

Each $x_i$ is the coefficient of a respective Fourier mode and the dynamics arise from physical processes as well as interactions between modes. The used parameter values are the same as in \cite{dorrington2023interaction} and may also be found in the provided code \cite{schoeller_2026_19711549}. Crucially, the system admits two fixed points, which were already identified in the original publication~\cite{charneyMultipleFlowEquilibria1979}. The associated flow fields vaguely resemble what is referred to as a ``blocking'' or a ``zonal'' atmospheric configuration, or weather regime, respectively. This is important because the system behaves considerably differently in the vicinity of the blocking fixed point compared to the zonal fixed point. Two dimensional projections of the attractor and its fixed points are shown in Fig.~\ref{fig:Hists}. 

\begin{figure}[H]
    \centering
    \includegraphics[width=0.95\textwidth]{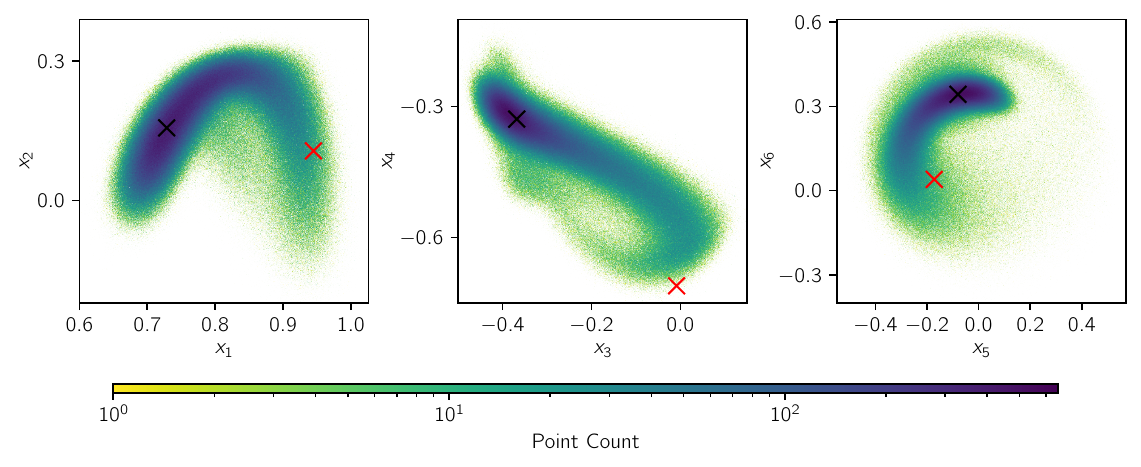}
    \caption{Same as Fig.~\ref{fig:Hists}, but for the system with additive noise of strength~$\sigma = 10^{-3}$. Note the different colour bar scaling to accommodate lower maximum bin count.}
    \label{fig:HistsStoch}
\end{figure}

We point out that the logarithmic colour scale in Fig.~\ref{fig:Hists} identifies very frequently visited orbits while the rest of the attractor is far less densely populated. The attractor's structure is typical for a system governed by intermittent chaos which can be understood to be emanating from a codimension-2 fold-Hopf bifurcation generating the strange attractor~\cite{deswartLoworderSpectralModels1988}.
% It should also be noted that the dimension estimate of the attractor is just below two; cf.~Fig.~\ref{fig:sigmaloglogcdv}.
As can be judged from the animations provided in the supplement, the dynamics are characterized by non-periodic oscillation between decelerated phases of highly predictable flow near the blocking fixed point and chaotic behavior everywhere else. %This is reflected also in the largest finite-time Lyapunov exponent (FTLE) as a measure of the local predictability of the system. The vicinity of the blocking fixed point is characterized by considerably lower maximum FTLEs compared to the rest of the attractor, as we show in Fig.~\ref{fig:FTLE100} in the appendix.

\paragraph{Stochastic perturbation and regime identification.}

The system is integrated in time using the explicit midpoint rule (second-order Runge--Kutta; RK2) along with the addition of a stochastic increment in a Euler--Maruyama style to incorporate noise:
\begin{equation}
\begin{aligned}
k_1 &= \Delta t^* \, f(z_n), \\
k_2 &= \Delta t^* \, f\Big(z_n + \frac{1}{2} k_1 \Big), \\
z_{n+1} &= z_n + k_2 + \sigma \sqrt{\Delta t^*}\, \xi_n,
\end{aligned}
\end{equation}
where $f(x)$ signifies the right-hand side of Eqn.~\eqref{eq:CdV}, and $\xi_n$ is a six-dimensional vector of iid normally distributed random variables that are white in time. All integrations were performed with a time step of $\Delta t^* = 10^{-5}$ to minimize numerical error. The trajectories analyzed were accordingly subsampled with the appropriate subsampling ratio to arrive at~$\Delta t$.  A fast \texttt{FORTRAN} implementation of this is provided in \cite{schoeller_2026_19711549}, which is a marginally modified version of the one provided by \cite{dorrington2023interaction}. A slightly slower, yet more convenient, parallelized, and just-in-time compiled implementation in \texttt{python} is also provided. 

We show a visualization of the attractor of the stochastic system for a selected noise strength $\sigma = 10^{-3}$ in Fig.~\ref{fig:HistsStoch}. The addition of noise blurs the detailed structure apparent in the deterministic case distributing the system's frequency of visits more evenly across the state space. This also decreases the histograms maximum bin count which is why we adjust the colour bar scaling compared to Fig.~\ref{fig:Hists}. More importantly, the noise also enables the system to visit regions that it was not able to reach in the deterministic case. This effect is especially pronounced in the vicinity of the blocking fixed point, where the system now resides far more frequently.

Equipped with the integrated trajectories of the CdV system, the next step in obtaining regime lifetimes is the definition of the regime. In \cite{dorrington2023interaction} the authors used a hidden Markov model to determine the blocking regime for each value of $\sigma$ separately, implying that the blocking regime changes with~$\sigma$. This was deemed necessary because of the changing shape of the attractor discussed above and the resulting set of blocking regimes across different noise levels exhibiting a pronounced shift in the mean position towards the blocking fixed point. However, since the Markov chain method introduced in Sec.~\ref{ssec:Markov_model} only needs a single deterministic trajectory and, therefore, only a single regime definition, we keep a single regime definition for all values of~$\sigma$.

To stay mathematically consistent, we refrain from defining the respective regime using a hidden Markov model, since our theoretical framework also allows us to objectively define a regime according to its \emph{dynamical properties}. Let $\mathcal{X}$ be a long integrated trajectory of the CdV model without noise. Based on $\mathcal{X}$, we construct the transition probability matrices $P_{ij}^{\sigma}$ for a small value of $\sigma$ as described in Sec.~\ref{ssec:Markov_model}. The regime $\mathcal{M} \subset \mathcal{X}$ is defined by a spectral clustering method with respect to the matrix $P^\sigma$: We use $P^\sigma$ to find almost-invariant sets under the action approximated by the transfer operator~\cite{dellnitz1999approximation}. This way, regimes are identified as subsets of the state space having a large dynamical stability. We calculate a small set of dominant left-eigenvectors of $P^{\sigma}$ (in our case four) and use them as embedding coordinates for our trajectory data. We then apply a k-means clustering on the embedded data to obtain a partition into clusters. Usually, the number of clusters would be determined according to a spectral gap, but since we are content with identifying a single blocking regime, we simply partition the data into two clusters.

\begin{figure}[]
    \centering
    \includegraphics[width=.6\textwidth]{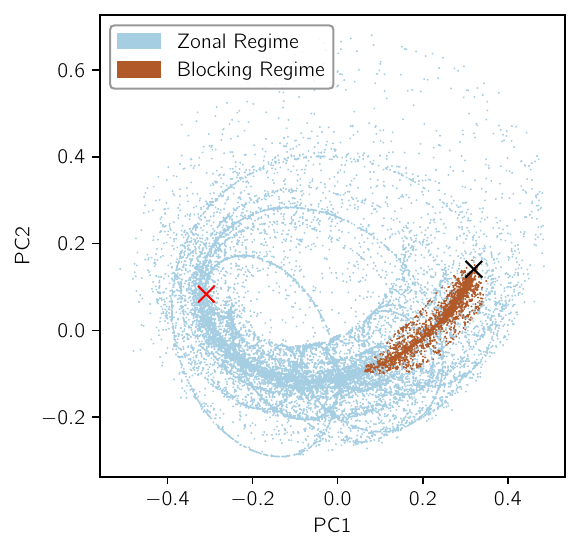}
    \caption{Assignment of $2\times 10^{4}$ points in a principal component projection from a deterministic trajectory with $\Delta t=10$ that has been partitioned according to the \texttt{k-means} clustering in the embedding explained in the text based on $\sigma = 10^{-2}$. Also indicated are the ``blocking'' fixed point (black) and the ``zonal'' fixed point (red).}
    \label{fig:Sets}
\end{figure}

The resulting clustering with $\sigma=10^{-2}$ is presented in Fig.~\ref{fig:Sets}. We choose this noise level since it lies in the range of values for $\sigma$ in which the Markov chain method is expected to work satisfyingly (see below) and since this is around the noise level where \cite{dorrington2023interaction} observed the maximum increase in regime lifetime. We will use the identified regime definition throughout the ensuing analysis. The algorithm identifies a cluster made up of a total of 4,989 points which we call the \emph{blocking regime}. This cluster is very similar to the blocking regime identified in \cite{dorrington2023interaction}, though less concentrated on the one-dimensional sub-manifold. Their hidden Markov model identifies 4,039 points in the blocking regime when applied to the same deterministic trajectory, of which 3,452 are shared with our regime.

\begin{figure}[]
    \centering
    \includegraphics[width=.6\textwidth]{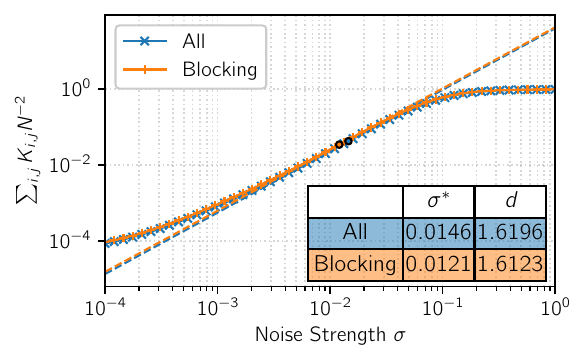}
    \caption{Normalized sum of affinities according to Eqn.~\eqref{eq:simiMat} plotted over the noise strength $\sigma$ for all points shown in Fig.~\ref{fig:Sets} (``All'') as well as for points either in the blocking regime or within a distance of $3\sqrt{2\sigma_{\max}^2 \Delta t}$ of any point in the blocking regime, where $\sigma_{\max}=0.03$ (``Blocking''). The dashed lines show the tangent at $\sigma = \sigma^{\ast}$ (indicated by open circles), the noise level maximizing the slope of the curve with the maximum slope~$d_\mathcal{X}$ (``All'') and $d_\mathcal{M}$ (``Blocking''). Both values are indicated for both sets in the table in the bottom right.}
    \label{fig:sigmaloglogcdv}
\end{figure}

To make sure that the sampled trajectory $\mathcal{X}$ is dense enough for our method to work, we again plot $\sum_{i,j}K^{\sigma}_{ij}$ against $\sigma$ in a doubly logarithmic plot in Fig.~\ref{fig:sigmaloglogcdv} and find a linear segment for the points relevant for the blocking regime in $2\times 10^{-3} \lesssim \sigma \lesssim 10^{-1}$. This is in satisfying agreement with the range of noise strengths we want to study and barely changes if instead all points in $\mathcal{X}$ are used. The maximum slope of roughly $d_\mathcal{M}=1.6$ suggests the blocking regime appears locally almost one-dimensionally, which confirms our intuition. 

\paragraph{Stochastic inertia.}

\begin{figure}[]
    \centering
    \includegraphics[width=.75\textwidth]{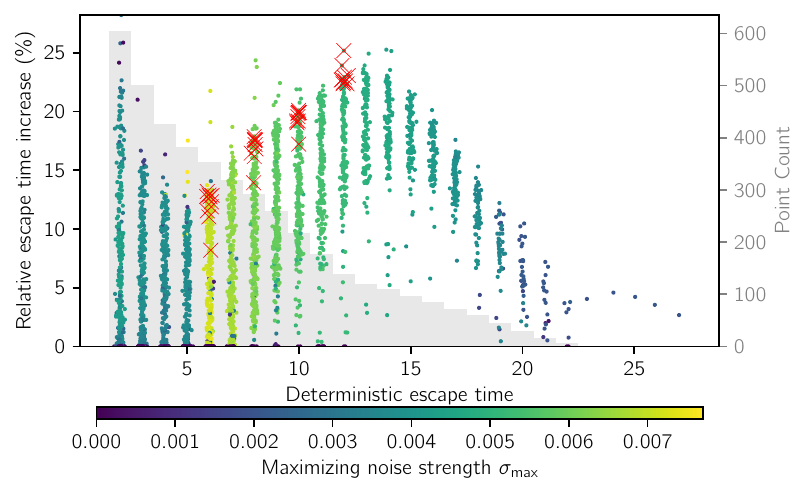}
    \caption{For each point in the blocking regime $\mathcal{M}$, its maximum relative increase in escape time (left y-axis) is shown against its deterministic escape time. A small random value is added to the x-location of each point to improve visibility. The points are coloured according to the noise strength $\sigma_{\mathrm{max}}$ at which the maximal escape time is achieved. Data for points within one deterministic step until escape are not shown. A grey histogram in the background shows the number of points that have a specific deterministic escape time (right y-axis). Indicated as red crosses are the points sampled for analysis shown in Fig.~\ref{fig:CdV_pointwise} due to their location on the attractor.}
    \label{fig:TOLifetimes_increase}
\end{figure}

Our goal is to verify the existence of stochastic inertia using our Markov chain method applied to the deterministic trajectory data $\mathcal{X}$. For $\sigma>0$ we use the transition probability matrix $P^\sigma$ to compute the expected escape times $\theta^\sigma$ of each point in $\mathcal{M}$ according to Eqn.~\eqref{eq:def_theta}. We will compare these expected escape times to those computed by a Monte Carlo simulation later and first focus in detail on the results of our Markov chain method.  Before showing the familiar curves of noise strength versus escape times for selected points, we stratify the results for all points according to their location in the blocking regime as well as their noise level of maximum escape time increase in Fig.~\ref{fig:TOLifetimes_increase}.

For each point $X_i \in \mathcal{M}$, we compute the expected escape times $\theta^\sigma_i$ for a range of values of~$\sigma$. The value of $\sigma$ that corresponds to the largest expected escape time is denoted by~$\sigma_\text{max}(X_i)$. The relative increase in the expected escape time is thus given by
\begin{equation}
    R(X_i) \coloneqq \frac{\theta_i^{\sigma_\text{max}(X_i)}}{\theta_i^0} - 1.
\end{equation}
Larger values of $R(X_i)$ correspond to stronger stochastic inertia. Fig.~\ref{fig:TOLifetimes_increase} shows the relative increase of the escape time for all points in $\mathcal{M}$. The x-axis displays the deterministic escape time $\theta_i^0$. Hence, points on the left side of the plot lie at the end of the blocking regime while points on the right of the plot lie at the entrance of the regime. In particular, the majority of points in the middle lie in the center of the blocking regime, which is where we expect the strongest effect of stochastic inertia. Indeed, Fig.~\ref{fig:TOLifetimes_increase} shows a large relative increase of the expected escape time for all points in the middle. This indicates that the Markov chain method is able to reliably predict the existence of stochastic inertia in the blocking regime of the CdV model. Note that we have deliberately clipped data from points within one deterministic step from escape, as results are very scattered here, both with respect to $R(X_i)$ and $\sigma_{\max}$.

Given that we are most confident in the results for points distant from entry and escape, we first compare escape times for such points with those obtained by the Monte Carlo method. We select the 10 points with the highest $R(X_i)$ and 12 deterministic steps away from escape and use their locations 12, 10, 8 and 6 steps before escape as initialization points in the following. They are indicated by red crosses in Fig.~\ref{fig:TOLifetimes_increase}.

\begin{figure}[]
    \centering
    \includegraphics[scale=.75]{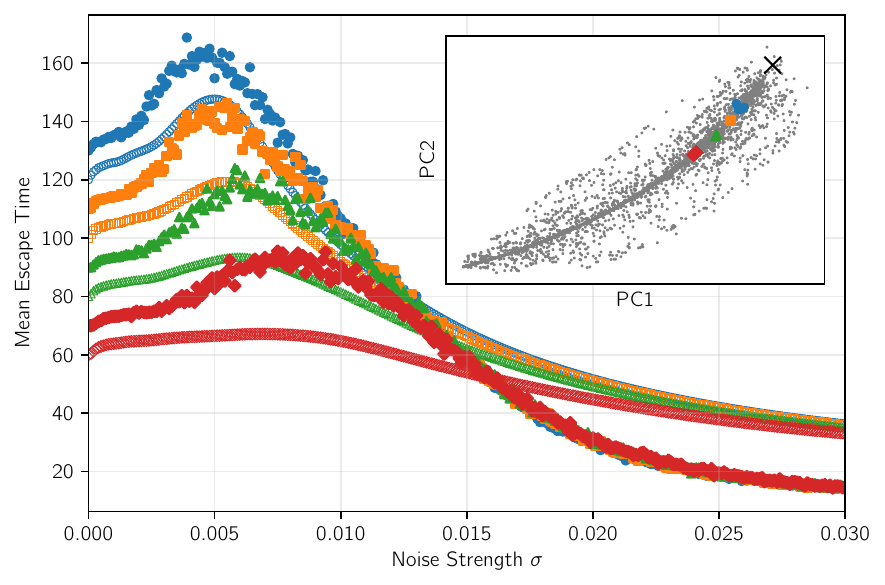}
    \caption{Pointwise escape times from the center of the blocking regime of the CdV system. The hollow markers belong to the expected escape times $\theta^{\sigma}_i$ calculated according to Eqn.~\eqref{eq:def_theta}. They are means over ten very close points and also indicated with red crosses in Fig.~\ref{fig:TOLifetimes_increase}. The solid markers belong to escape times calculated in a Monte Carlo fashion for the same sets of points with $100$ realizations for each point and noise strength. The location of the points are also shown in the inset along with all the other points in $\mathcal{M}$ in the same principal component projection used in Fig.~\ref{fig:Sets} and the ``blocking'' fixed point as a black cross.}
    \label{fig:CdV_pointwise}
\end{figure}

For all of these points that lie well within the blocking regime, the average escape time computed using the Markov chain method matches well with the average escape time estimated by Monte Carlo simulations and clearly exhibits stochastic inertia; see Fig.~\ref{fig:CdV_pointwise}. The difference between the Markov chain method and the Monte Carlo estimations is probably due to the Markov chain method's restriction to the deterministic attractor, as discussed in~Sec.~\ref{sec:lim} f). More specifically, we consider a Monte Carlo trajectory to have escaped the regime once there is no point from $\mathcal{M}$ within a distance of $10^{-1}$. This also explains the difference in escape time in the deterministic case $\sigma=0$; it is simply the time needed for the deterministic trajectory to travel a distance of $10^{-1}$. 

\begin{figure}[]
    \centering
    \includegraphics[scale=.75]{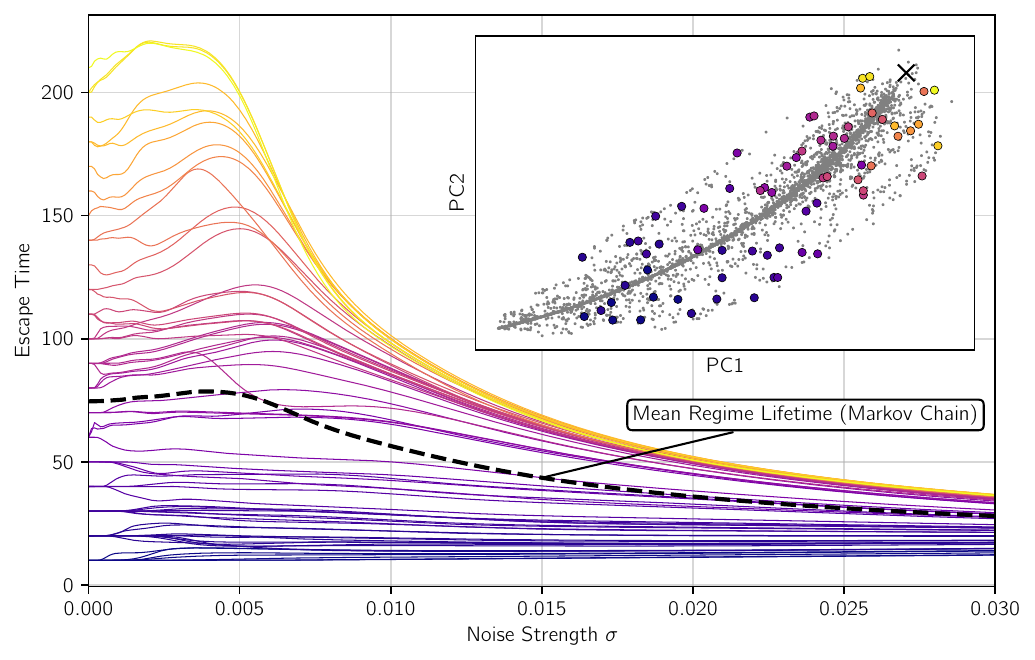}
    \caption{Expected escape times $\theta^{\sigma}$ across noise strength $\sigma$ calculated according to Eqn.~\eqref{eq:def_theta} for points of the deterministic trajectory that have just entered the blocking regime $\mathcal{M}$. Only every tenth entry point is shown for better visibility. The curves are coloured according to their deterministic escape time ($\sigma=0$). The location of these entry points are also shown in the inset along with all the other points in $\mathcal{M}$ in the same principal component projection used in Fig.~\ref{fig:Sets} and the ``blocking'' fixed point as a black cross. The black dashed line shows the mean over all curves, which is an approximation of the expected regime life time $E^{\sigma}$.}
    \label{fig:CdV_LifetimesEntries}
\end{figure}

Given the complexity of the CdV system, the agreement between the two methods is satisfying, but the very dense, almost one-dimensional submanifold we focus on in Fig.~\ref{fig:CdV_pointwise} is the region where the method is most likely to work. Moreover, the set of initialization points used differ substantially from the actual entry distribution. By contrast, the escape time for different noise levels for every tenth deterministic entry points (subsampled for visibility purposes) is shown in Fig.~\ref{fig:CdV_LifetimesEntries}. Comparing the inset of that figure to the one of Fig.~\ref{fig:CdV_pointwise} proves that the entry points' locations are spread out across the blocking regime. The colour of the points and the associated curves in Fig.~\ref{fig:CdV_LifetimesEntries} indicate points travel through the regime away from the fixed point. Indeed, the convergence of the individual curves also proves that points in the regime ultimately reach the one-dimensional submanifold in its center making their location and, hence, escape time virtually indistinguishable. 

Fig.~\ref{fig:CdV_LifetimesEntries}, in accordance with Fig.~\ref{fig:TOLifetimes_increase}, demonstrates that stochastic inertia is found by our method for all entry points that are sufficiently far away from (deterministic) escape. Following the method described in Sec.~\ref{ssec:Markov_model}, we need to estimate the entry distribution $\nu^\sigma$ as defined in Eqn.~\eqref{eq:def_nu} in order to compute the expected regime lifetimes. As described in Sec.~\ref{sec:lim} d), a crucial problem in the computation of $\nu^\sigma$ is that the computational complexity scales cubically in the size $|\mathcal{X}|$ of the entire point cloud. Meanwhile, computing the escape times $\theta^\sigma$ only scales cubically in the size $|\mathcal{M}|$ of the blocking regime. When the resolution of $\mathcal{M}$ is sufficiently large for our method to work, $|\mathcal{X}|$ already becomes too large to compute $\nu^\sigma$ by conventional methods. Since this is not the focus of our work, we do not discuss potential workarounds. As we discuss below, estimating the expected regime lifetimes using the Markov chain method faces other problems than just computational complexity. 
To have a value of $E^\sigma$ that can be used for comparison, we choose $\nu^\sigma$ to be the uniform distribution over all entry points of the deterministic trajectory. For small $\sigma$ this is a reasonable assumption since $\nu^0$ is exactly the uniform distribution over all entry points. The resulting expected regime lifetimes $E^\sigma$ are displayed in Fig.~\ref{fig:CdV_LifetimesEntries} as a black dashed line.

So far, we have shown that our proposed Markov chain method is able to reproduce escape time estimates obtained with the Monte Carlo method for points in the center of the blocking regime and that it can predict stochastic inertia in the regime lifetime of the blocking regime. The final test left is to compare the regime lifetimes of the blocking regime between the Markov chain and the Monte Carlo methods. As we will show, this is a challenging task given the different parts of the state space the two methods ``see'' and their respective computational limitations.

Given the definition of the blocking regime $\mathcal{M}$ shown in Fig.~\ref{fig:Sets}, we want to compute the average regime lifetimes using Monte Carlo simulations for different values of the noise strength~$\sigma$. Since $\mathcal{M}$ is only a finite set of points in $\R^6$ and not an open set, it is non-trivial to decide whether a point $z_n^\sigma$ of the stochastic Monte Carlo simulation lies inside or outside of the blocking regime. We say that a point $z_n^\sigma$ lies inside of the blocking regime, if and only if at least half of its 10 nearest neighbors within $\mathcal{X}$ belong to~$\mathcal{M}$. We use \texttt{scipy}s \texttt{cKDTree} for this \cite{2020SciPy-NMeth}. This regime assignment is slightly different from the one employed for the comparison shown in Fig.~\ref{fig:CdV_pointwise}, but considerably more convenient computationally because only ten neighbors have to be found instead of all of them. The results of both methods always agreed in our tests (not shown). 

\begin{figure}[]
    \centering
    \includegraphics[width=.75\textwidth]{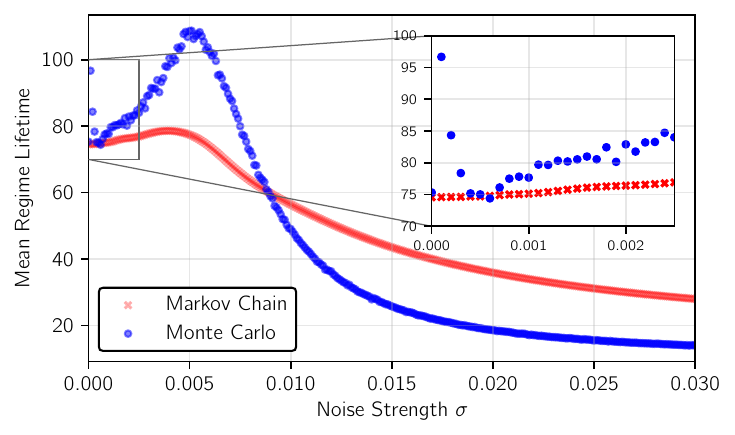}
    \caption{Mean regime life times for each noise strength $\sigma$ obtained from either the Monte Carlo method and the Markov chain approach. The Monte Carlo trajectory is computed with a step size of $10^{-5}$. We subsample the trajectories such that both methods have a time step of $\Delta t = 10$, but the Markov chain method uses $2 \times 10^4$ steps compared to $2 \times 10^6$ steps for the Monte Carlo method. The Monte Carlo trajectories feature roughly $10^4$ entries with $669$ entries  for the Markov Chain method (every tenth of which depicted in Fig.~\ref{fig:CdV_LifetimesEntries}). The inset details the low-noise behavior of the two.}
    \label{fig:CdVLifetInc}
\end{figure}

The average lifetimes of these Monte Carlo simulations for a range of values of $\sigma$ are displayed in Fig.~\ref{fig:CdVLifetInc}. Both methods indicate stochastic inertia at a similar noise level, though with considerably differing intensity. A peculiar intermittent increase is also observed for some of the very low noise cases for the Monte Carlo method. We expect this to be a numerical artifact rather than a meaningful effect. We discuss potential causes for the mismatch between the two curves. 

The blocking regime, as shown in Fig.~\ref{fig:Sets} is relatively small, such that sampling long regime lifetimes with a Monte Carlo method constitutes sampling rare events. This effect is strongly amplified by the fact that the shape of the attractor changes under the addition of even small noise; compare Fig.~\ref{fig:HistsStoch}. The Monte Carlo simulation of Dorrington and Palmer \cite[Figure 6]{dorrington2023interaction} differs from our simulations only by the choice of time step (they used $\Delta t=1$) and the definition of the blocking regime. Their experiments indicate the existence of stochastic inertia extends up to higher noise levels and increases in magnitude in that range. We note that their approach reveals a similar volatility at very low noise levels and that the mean regime lifetime agrees well with our results for noise strength up to approximately $\sigma = 0.005$ -- the maximum in our case. In contrast to the mean regime lifetimes calculated using our Markov chain framework, we have observed that the regime lifetimes of the Monte Carlo simulation are quite sensitive with respect to the step size $\Delta t$. We emphasize that to compute regime lifetimes, it is not always desirable to choose $\Delta t$ as small as possible to increase accuracy. As $\Delta t$ goes to zero, the stochastic trajectory converges to the trajectory of an SDE which does not have well-defined average regime lifetimes; see the discussion in Section \ref{sec:regime_lifetimes}. In addition, the choice of $\Delta t$ also influences the geometry of the blocking regime and, accordingly, the regime lifetimes calculated using the Monte Carlo method. As an example, we performed our analysis for $\Delta t = 1$. The corresponding results, analogues of Figs.~\ref{fig:Sets} and~\ref{fig:CdVLifetInc} are included in App.~\ref{app:Figs} as Fig.~\ref{fig:SetsDT1} and \ref{fig:CdVLifetIncDT1}. In Sec.~\ref{sec:lim} c) we argued that the time step $\Delta t$ should not bee too large. However, comparing Figs.~\ref{fig:Sets} and~\ref{fig:CdVLifetInc}, to Figs.~\ref{fig:SetsDT1} and \ref{fig:CdVLifetIncDT1}, it seems that the results of the Markov chain method are more robust to changes in $\Delta t$ than those of the Monte Carlo method.

A trajectory with smaller time steps will lead to a set $\mathcal{M}$ with fewer passages through it, implying a more skeletal structure. With our method of assigning points in noisy trajectories from the Monte Carlo scheme to each of the two regimes defined by the deterministic trajectory using a nearest-neighbour classification, this more skeletal structure leads to an increase in the number of (and, in turn, shorter) visits to the regime, biasing the mean lifetime downwards. This effect becomes more important for stronger stochastic forcing. The hidden Markov model regime definition employed by Dorrington and Palmer, on the other hand, prevents such a behavior by redefining the regimes for every trajectory. This trade-off arises as a consequence of the more fundamental decision of what exactly to consider the regime and it  relates to another vital discrepancy between the two procedures.

The Markov chain framework uses the state space $\mathcal{X}$ for all values of~$\sigma$. As mentioned before, and shown in Fig.~\ref{fig:HistsStoch}, even for very small noise the stochastic dynamics explores regions of the state space that lie outside of the deterministic attractor, and the stationary distribution changes significantly. Hence, we are confronted with the issues raised in Sec.~\ref{sec:lim} e) and f), and the entry distribution for intermediate $\sigma$ cannot be captured accurately by a probability distribution $\nu^\sigma $ on $\mathcal{X}$. The Monte Carlo method does not face this limitation. 

Finally, the density of the blocking regime $\mathcal{M}$ is highest in its center and lowest at its boundary. Our analysis revealed that , as a consequence, our Markov chain approach seems to perform worse at the boundary of the blocking regime than in the center. Since the escape times of the boundary points are exactly what is used for the regime lifetimes $E^\sigma$, the results lose accuracy.

\section{Conclusion}\label{sec:fin}

We proposed a method to numerically investigate noise-induced enhancement of regime lifetimes (which we called \emph{stochastic inertia}) without simulating the underlying system for different noise levels. Given a ``base'' trajectory with its own noise level (usually, but not necessarily, zero), we emulate additional noise by defining a suitable Markov chain on the points given by this base trajectory. This way, we were able to investigate the system's response to multiple noise levels without additional information about or simulation of the system. By testing our method against a detailed Monte Carlo simulation we delivered numerical evidence that a one-dimensional unstable fixed point exhibits stochastic inertia, while a three-dimensional one does not. Theoretical justification of these observations will appear elsewhere.

Additionally, we investigated a paradigmatic model for atmospheric dynamics; the Charney--deVore system. Our analysis deviated from the previous study~\cite{dorrington2023interaction} considering the phenomenon of regime-lifetime extension under increasing noise, since we defined the blocking regime by its strong dynamical persistence (almost-invariance). Yet, the effect of stochastic inertia was still observable.

For the one- and three-dimensional toy examples, our Markov chain method was able to accurately predict the average escape times and expected regime lifetimes. For the Charney--deVore system, our method was able to detect stochastic inertia and accurately compute the escape times of points in the center of the blocking regime, however, the method struggled to accurately compute the expected regime lifetimes.
% Additionally, our analysis seems to be sensitive to changes in the step size $\Delta t$. 
Comparing the results of Sec.~\ref{sec:CdV} to the corresponding figures in App.~\ref{app:Figs} indicates that the Monte Carlo simulation is sensitive to changes in $\Delta t$ while our Markov chain model seems less sensitive to such changes. Future efforts could target current limitations, such as potential changes in the attractor and the stationary distribution of the system due to changing noise levels, which are not accounted for in our method. An additional direction is to better understand how parameter and design choices impact results and, in particular, differences between the two methods presented. 

\subsection*{Author Contributions}

All authors contributed to the conceptual
framework of the paper. RC developed the theory and was the main author of Sect.~\ref{sec:meth}. HS carried out the computations, was responsible for creating the figures and was the main author of Sect.~\ref{sec:app}. All authors contributed to revision of the manuscript.

\subsection*{Acknowledgements}
This research has been supported by Deutsche Forschungsgemeinschaft (DFG) through grant CRC 1114 Scaling Cascades in Complex Systems, Project Number 235221301, Project A08, “Characterization and Prediction of Quasi-Stationary Atmospheric States”.
ME additionally thanks the Dutch Research Council (NWO) for support under the grant VI.Vidi.233.133.

\appendix

\section{Additional Figures}\label{app:Figs}

\begin{figure}[H]
    \centering
    \includegraphics[width=.6\textwidth]{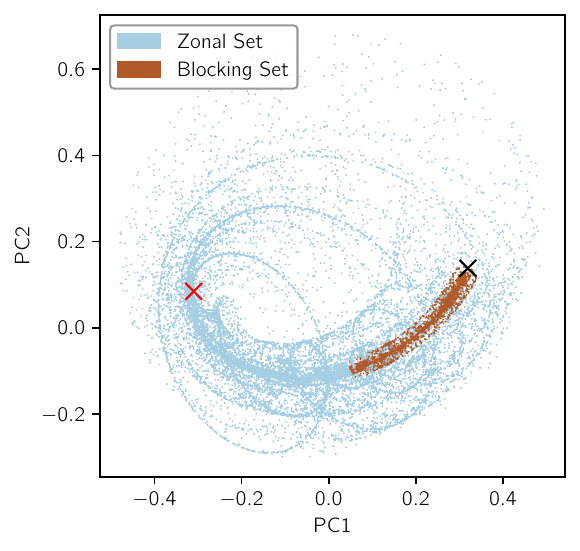}
    \caption{Same as Fig.~\ref{fig:Sets} but using a trajectory with $\Delta t=1$.}
    \label{fig:SetsDT1}
\end{figure}

\begin{figure}[H]
    \centering
    \includegraphics[width=.75\textwidth]{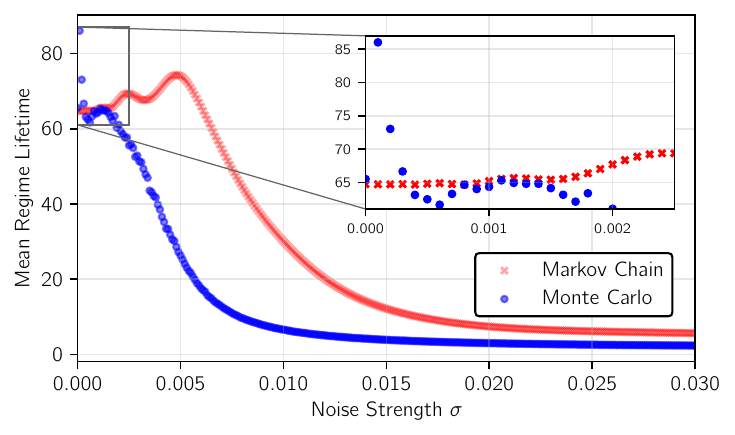}
    \caption{Same as Fig.~\ref{fig:CdVLifetInc}, but for trajectories with $\Delta t = 1$.
    }
    \label{fig:CdVLifetIncDT1}
\end{figure}

\section{Remarks on the evaluation of eqn.~(\ref{eq:analytic_1d})}
\label{app:NumRem}
Starting from \eqref{eq:analytic_1d}, we rewrite the exponent
\[
-\frac{y^2 - z^2}{\sigma^2} = -\frac{y^2}{\sigma^2} + \frac{z^2}{\sigma^2},
\]
which yields
\[
\E[\tau^\sigma (x)] 
= \frac{\displaystyle \int_{-1}^x e^{-\frac{2y^2}{\sigma^2}}
\int_{-1}^y e^{\frac{z^2}{\sigma^2}}\,dz\,dy}
{\displaystyle \frac{1}{2}\sigma^2 \int_{-1}^1 e^{-\frac{y^2}{\sigma^2}}\,dy}.
\]

Introducing the rescaling \(u = y/\sigma\), \(w = z/\sigma\), the inner integral becomes
\[
\int e^{w^2}\,dw = \frac{\sqrt{\pi}}{2}\,\mathrm{erfi}(w),
\]
so that, after inserting limits and simplifying constants with the denominator, the expression reduces to a single integral
\[
\E[\tau^\sigma(x)] = \int_{u_1}^{u_2} g(u)\,du,
\qquad 
g(u) = -\sqrt{\pi}\,e^{-u^2}\,\mathrm{erfi}(u),
\]
with
\[
u_1 = -\frac{1}{\sqrt{2}\,\sigma}, \qquad
u_2 = \frac{x}{\sqrt{2}\,\sigma}.
\]

For numerical evaluation, this integral is computed using adaptive quadrature. To ensure stability, the domain is split into a central region \(|u|\le U\), where \(g(u)\) is evaluated directly, and tail regions \(|u|>U\), where the asymptotic form \(g(u)\sim -1/u\) is used. The contributions are then summed.

\bibliographystyle{myalpha}
\bibliography{bibliography.bib}

\end{document}